\documentclass{article}

\usepackage[english]{babel}

\usepackage{amsmath}
\usepackage{amsfonts}
\usepackage{tikz}
\usepackage{enumitem}
\usepackage{placeins}

\def\vec#1{\ensuremath{\mathchoice
    {\mbox{\boldmath$\displaystyle\mathbf{#1}$}}
    {\mbox{\boldmath$\textstyle\mathbf{#1}$}}
    {\mbox{\boldmath$\scriptstyle\mathbf{#1}$}}
    {\mbox{\boldmath$\scriptscriptstyle\mathbf{#1}$}}}}%

\newcommand{\vect}[1]{\vec{#1}}
\def\tens#1{\relax\ifmmode\mathsf{#1}\else\textsf{#1}\fi}
\newcommand{\matr}[1]{\tens{#1}}
\newcommand{\transpose}[1]{{#1}^{\text{T}}}

\def\titre#1{\begin{center}{\Large{\bf #1}}\end{center}}

\def\motscles#1{%
	\ifx#1\IsUndefined\relax\else\noindent{\normalsize{\bf Keywords:}} #1\\ \fi}

\begin{document}

\titre{A conservative coupling algorithm between a
  compressible flow and a rigid body using an Embedded Boundary method}

\begin{center}{\bf Laurent~Monasse$^{\textrm{1, 2, 3}}$, Virginie~Daru$^{\textrm{3, 4}}$, Christian~Mariotti$^{\textrm{1}}$, Serge~Piperno$^{\textrm{2}}$, Christian~Tenaud$^{\textrm{3}}$} \end{center}
\begin{center}
  $^{\textrm{1}}$ CEA-DAM-DIF 91297 Arpajon, France \\
  email: {\tt \{laurent.monasse, christian.mariotti\}@cea.fr}
\end{center}
\begin{center}
  $^{\textrm{2}}$ Universit\'e Paris-Est, CERMICS (ENPC), \\
  77455 Marne la Vall\'ee cedex, France \\
  email: {\tt \{monassel, piperno\}@cermics.enpc.fr}
\end{center}
\begin{center}
  $^{\textrm{3}}$ LIMSI-CNRS, 91403 Orsay, France \\
  email: {\tt \{laurent.monasse, virginie.daru, christian.tenaud\}@limsi.fr}
\end{center}
\begin{center}
  $^{\textrm{4}}$ Lab. DynFluid, Ensam, 75013 Paris, France \\
  email: {\tt virginie.daru@ensam.eu}
\end{center}

\rule{\linewidth}{.5pt}
\\

  \textbf{ABSTRACT}

This paper deals with a new solid-fluid coupling algorithm
between a rigid body and an unsteady compressible fluid flow, using an
Embedded Boundary method. The
coupling with a rigid body is a first step towards the coupling with
a Discrete Element method. The flow is computed using a Finite
Volume approach on a Cartesian grid. The expression of numerical
fluxes does not affect the general coupling algorithm and we use a
one-step high-order scheme proposed by Daru and Tenaud [Daru
V,Tenaud C., J. Comput. Phys. 2004]. The Embedded Boundary method is
used to integrate the presence of a solid boundary in the fluid. The
coupling algorithm is totally explicit and ensures exact mass
conservation and a balance of momentum and energy between the fluid
and the solid. It is shown that the scheme preserves uniform
movement of both fluid and solid and introduces no numerical
boundary roughness. The efficiency of the method is demonstrated on
challenging one- and two-dimensional benchmarks.

\\
\rule{\linewidth}{.5pt}
\parskip = 0.mm
\parindent = 5mm

\renewcommand{\labelitemi}{$-$}

\section{Introduction}

This work is devoted to the development  of a coupling method for
fluid-structure interaction in the compressible case. We intend to
simulate transient dynamics problems, such as the impact of shock
waves onto a structure, with possible fracturing causing the
ultimate breaking of the structure. An inviscid fluid flow model is
considered, being convenient for treating such short time scale
phenomena. The simulation of fluid-structure interaction problems is
often computationally challenging due to the generally different
numerical methods used for solids and fluids and the instability
that may occur when coupling these methods. Monolithic methods have
been employed, using an Eulerian formulation for both the solid and
the fluid (for instance, the diffusive interface method
\cite{Gavrilyuk,Abgrall}), or a Lagrangian formulation for both the
fluid and the solid (for example, the PFEM method \cite{ref8}), but
in general, most solid solvers use Lagrangian formulations and fluid
solvers use Eulerian formulations. In this paper we consider the
coupling of a Lagrangian solid solver with an Eulerian fluid solver.

For the coupling in space, a possible choice is to
deform the fluid domain in order to follow 
the movement of the solid boundary: the Arbitrary
Lagrangian-Eulerian (ALE) method has been developed to that end. It
has been widely used for incompressible \cite{Halleux,Fernandez} and
compressible \cite{Farhat} fluid-structure interaction. However,
when solid impact or fracture occur, ALE methods are faced with a
change of topology in the fluid domain that requires remeshing and
projection of the fluid state on the new mesh, which are costly
and error prone procedures. Moreover remeshing is poorly adapted to
load balancing for parallel computations.

In order to allow for easier fracturing of the solid, we instead
choose a method based on fictitious domains that solves the fluid
flow on a fixed Eulerian mesh, on which a Lagrangian solid body is
superimposed. A special treatment is then applied on fluid cells
near the boundary and inside the solid. Different types of
fictitious domain methods have been developed over the last thirty
years. They can roughly be classified in three main classes:
penalization methods, interpolation methods and conservative
methods. Among penalization methods, the Immersed Boundary method is
certainly the best known and most widely used for fluid-structure
interaction. It was originally introduced by Peskin for
incompressible blood flows \cite{Peskin,IBM16}. The solid boundaries
deform under the action of the fluid velocity, and the presence of
the solid adds forces in the fluid formulation that enforce the
impermeability of the solid. However, Xu and Wang have pointed out
some numerical leaking of fluid into the solid \cite{IBM4}.
Following Leveque and Li \cite{Leveque,Leveque2}, they advocate the
use of the Immersed Interface method, which incorporates jump
conditions in the finite differences used. However, the absence of
fluid mass loss is still not ensured exactly. In a different
approach, Olovsson {\it et al.} \cite{Olovsson,Aquelet} couple an
Eulerian and a Lagrangian method by penalizing the penetration of
the solid into the fluid by a damped spring force. As the stiffness
of the spring goes to infinity, the penetration goes to zero. Boiron
{\it et al.} \cite{ref31} and Paccou {\it et al.} \cite{ref32}
consider the solid as a porous medium, using a Brinkman porosity
model. As the porosity goes to zero, the solid becomes impermeable.
However, in both cases, as the stiffness grows or the porosity
decreases, the use of implicit schemes is mandatory to avoid the
severe stability condition of explicit schemes \cite{ref31,ref32}.
For the high speed phenomena we consider, we use explicit solid and
fluid solvers and an explicit coupling algorithm is better suited in
order to avoid costly iterative procedures.

A second class of fictitious domain methods consists in enforcing the
boundary conditions through interpolations in the vicinity of the
boundary, using the exact values taken by the fluid on the boundary
\cite{MohdYusof,Fadlun}. The method seems to be very versatile,
being used with incompressible Navier-Stokes
\cite{MohdYusof,Fadlun}, Reynolds-averaged Navier-Stokes
\cite{IBM3,IBM7,IBM5}, turbulent boundary layer laws \cite{IBM1} and
compressible Navier-Stokes \cite{IBM6}. The Ghost Fluid method
developed by Fedkiw {\it et al.} \cite{IBM25,IBM31} relies on the
same type of principle for compressible fluids. The interface is
tracked using a level-set function, and conditions are applied on
both sides of the interface to interpolate the boundary conditions.
The advantage of these methods is that they do not suffer from
additional time-step restriction due to stability, and the order of
accuracy of the boundary conditions can be set {\it a priori}.
However, the interpolation does not ensure the conservation of mass,
momentum and energy in the system. This can cause problems when dealing
with shock waves interacting with solids.

In this article, we rather consider the third class of
  conservative fictitious domain methods, which seems the most adequate
  framework to develop our coupling algorithm. These methods are
  generally referred to as Embedded Boundary methods, and they rely on a
  modified integration of the numerical fluid fluxes in the cells cut by
  the solid boundary \cite{IBM18,IBM30,IBM28,IBM10}. The
original idea of the method can be traced back to Noh's CEL code
\cite{Noh}. The new contribution of the present work
  consists in the coupling algorithm and its properties. The
  Embedded Boundary method that we use is essentially identical to
  previous works \cite{IBM18,IBM30,IBM28,IBM10}. The different versions
  of the Embedded Boundary method mainly differ in the way the
  stability condition arising from small cut-cells is
  enforced and we develop here a slightly different procedure in order
  to deal with solid boundaries coming close to each other. The method
can be 
implemented independently from the
time integration scheme used for the fluid, whether based on
space-time splitting or multi-level time integration. Conservative
fictitious domain methods have proven to give satisfactory
conservation results for inviscid compressible flows in the case of
static solid boundaries. Nevertheless, to our knowledge,
conservation issues of the coupling have not been studied in the
case of moving solids. We establish new conservation results in such
a case. Our coupling method is designed to be capable of treating
the general case of moving deformable bodies. In the present work,
however, we only consider non-deformable (rigid) solid bodies. The
case of deformable bodies is the object of ongoing work.

The fluid and solid solvers that we consider were chosen according
to their ability to deal with shock waves and fracturing solids. The
solid solver is based on a Discrete Element method, implemented in a
code named Mka3D in the CEA \cite{Mariotti}.  It can handle
elasticity as well as fracture and impact of solids. Solids are
discretized into polyhedral particles, which interact through
well-designed forces and torques. The particles have a rigid-body
motion, and fracture is treated in a straightforward way by removing
the physical cohesion between particles. The work reported in this
article is a first step towards the coupling with the Mka3D code.
The time integration scheme used by Mka3D (Verlet for displacement
of the center of mass and RATTLE for rotation \cite{Monasse}) is
retained for the rigid body treatment. Concerning the fluid solver,
we use a Cartesian grid explicit finite volume method, based on the
high-order one-step monotonicity-preserving scheme developed in
\cite{Daru2004} and space  time splitting. However we emphasize
that our coupling method is independent from both the Discrete
Element method (as long as a solid interface is defined) and the
numerical scheme used for the fluid calculation.

The article is organized as follows: we first present briefly the
solid and fluid methods in section \ref{sec:solidfluid}. In sections
\ref{sec:cutcells} and \ref{sec:description}, we describe the
proposed explicit coupling procedure between the fluid and the
moving solid in the framework of an Embedded Boundary method. The
analysis of the conservation properties of the coupling is reported
in section \ref{sec:conservation}, where we show that mass, momentum
and energy of the solid-fluid system are exactly preserved.  In
section \ref{sec:constant}, we demonstrate results about the
preservation on a discrete level of two solid-fluid systems in
uniform movement. Finally, we illustrate the efficiency and accuracy
of the method on one and two-dimensional static and dynamic
benchmarks in section \ref{sec:numeric}.

\section{Solid and fluid discretization methods}
\label{sec:solidfluid}

\subsection{Solid time-discretization method}
\label{sec:solide}

We consider a non-deformable solid (rigid body). The position and
velocity of the solid are given, respectively, by the position of
its center of mass $\vect{X}$, the rotation matrix $\matr{Q}$, the
velocity of the center of mass $\vect{V}$ and the angular momentum
matrix $\matr{P}$. The physical characteristics of the solid are its
mass $m$ and its matrix of inertia $\matr{R}$ which, in the inertial
frame, is a diagonal matrix with the principal moments of inertia
$I_1$, $I_2$ and $I_3$ on the diagonal. Here, we instead use the
diagonal matrix $\matr{D}=\text{diag}(d_1,d_2,d_3)$, where:
\begin{equation*}
  \forall i\in\{1,2,3\}, d_i=\frac{I_1+I_2+I_3}{2}-I_i.
\end{equation*}
The angular momentum matrix $\matr{P}$ can be related to the usual
angular velocity vector $\vect{\Omega}$ by the relation $\matr{P} =
\matr{D}\matr{j}(\vect{\Omega})\matr{Q}$, where  the map $\matr{j} :
\mathbb{R}^3 \rightarrow \mathbb{R}^{3\times3}$ is defined such
that:
\begin{equation*}
  \forall \vect{x}\in\mathbb{R}^3,\; \forall
  \vect{y}\in\mathbb{R}^3,\;
  \matr{j}(\vect{x})\cdot\vect{y}=\vect{x}\wedge\vect{y}.
\end{equation*}

Let us denote by $\vect{F}$ and $\vect{\mathcal{M}}$ the external
forces and torques acting on the solid, and by $\Delta t$ the
time-step. In order to preserve the energy of the solid over time
integration of the movement, we choose a symplectic second-order
scheme for constrained Hamiltonian systems, the RATTLE scheme
\cite{Hairer}:
\begin{align}
   \vect{V}^{n+\frac{1}{2}} &= \vect{V}^{n} + \frac{\Delta
     t}{2m}\vect{F}^n, \label{eqn:vitesse1} \\
   \vect{X}^{n+1} &= \vect{X}^n + \Delta t\vect{V}^{n+\frac{1}{2}}, \label{eqn:position}\\
   \matr{P}^{n+\frac{1}{2}} &= \matr{P}^{n} + \frac{\Delta t}{4}\matr{j}(\vect{\mathcal{M}}^n)\matr{Q}^n + \frac{\Delta t}{2}\matr{\Lambda}^n\matr{Q}^n, \label{eqn:vitesserot1}\\
   \matr{Q}^{n+1} &= \matr{Q}^n + \Delta t\matr{P}^{n+\frac{1}{2}}\matr{D}^{-1}, \label{eqn:rotation}
 \end{align}
 \begin{equation}
   \text{with }\matr{\Lambda}^n\text{ such that } \transpose{(\matr{Q}^{n+1})}\matr{Q}^{n+1}=\matr{I}, \label{eqn:contrainte1}
 \end{equation}
 \begin{align}
   \vect{V}^{n+1} &= \vect{V}^{n+\frac{1}{2}} + \frac{\Delta t}{2m}\vect{F}^{n+1}, \label{eqn:vitesse2}\\
   \matr{P}^{n+1} &= \matr{P}^{n+\frac{1}{2}} + \frac{\Delta t}{4}\matr{j}(\vect{\mathcal{M}}^{n+1})\matr{Q}^{n+1} + \frac{\Delta t}{2}\matr{\tilde{\Lambda}}^{n+1}\matr{Q}^{n+1}, \label{eqn:vitesserot2}
 \end{align}
 \begin{equation}
   \text{with }\matr{\tilde{\Lambda}}^{n+1}\text{ such that } \transpose{(\matr{Q}^{n+1})}\matr{P}^{n+1}\matr{D}^{-1}+\matr{D}^{-1}\transpose{(\matr{P}^{n+1})}\matr{Q}^{n+1}=\matr{0}. \label{eqn:contrainte2}
 \end{equation}
The symmetric matrices $\matr{\Lambda}^n$ and
$\matr{\tilde{\Lambda}}^{n+1}$ play the role of Lagrange multipliers
for the constraints on matrices $\matr{Q}^{n+1}$ and
$\matr{P}^{n+1}$.

The scheme makes use of the velocity at half time-step
$\vect{V}^{n+\frac{1}{2}}$, which is constant during the time-step.
Let us now consider the angular velocity. For a rigid solid, we have
for all points  $\vect{x}$:
\begin{equation*}
  \vect{X}-\vect{x} =
  \matr{Q}\cdot\left(\vect{X}^0-\vect{x}^0\right),
\end{equation*}
$\vect{X}^0$ and $\vect{x}^0$ being material points of the solid at
initial time. Using the identity
$\vect{\Omega}\wedge(\matr{Q}\vect{x})=\matr{P}\matr{D}^{-1}\vect{x}$
for all $\vect{x}$, the velocity at point $\vect{x}$ can be written
as:
\begin{equation*}
  \vect{V}(\vect{x}) = \vect{V} +
  \matr{P}\matr{D}^{-1}\cdot\left(\vect{X}^0-\vect{x}^0\right)
\end{equation*}
which is more convenient for use in the time scheme. In analogy with
displacement, we consider $\matr{P}^{n+\frac{1}{2}}$ as constant
during the time-step, and we define the velocity of point $\vect{x}$
at half time-step $(n+\frac{1}{2})\Delta t$:
\begin{equation*}
  \vect{V}^{n+\frac{1}{2}}(\vect{x}) = \vect{V}^{n+\frac{1}{2}} +
  \matr{P}^{n+\frac{1}{2}}\matr{D}^{-1}\cdot\left(\vect{X}^0-\vect{x}^0\right).
\end{equation*}

\subsection{Fluid discretization method}

The problem of the interaction of shock waves with solid surfaces
can be at first studied using an inviscid fluid model. In this work,
we consider inviscid compressible flows, which follow the Euler
equations:
\begin{equation*}
  \vect{w}_t+\nabla\cdot\vect{f}(\vect{w}) = \vect{0},
\end{equation*}
where $\vect{w}=\transpose{(\rho,\rho \vect{u},\rho E)}$ is the vector
of the conservative variables, and $\vect{f}(\vect{w})$ is the Euler
flux:
\begin{equation*}
  \vect{f} = \left(
    \begin{array}{c}
      \rho\vect{u} \\
      \rho\vect{u}\otimes\vect{u}+p\matr{I}\\
      (\rho E+p)\vect{u}
    \end{array} \right),
\end{equation*}
where the pressure $p$ is given by a perfect gas law: $p = (\gamma-1)\left(\rho E-\frac{1}{2}\rho\vect{u}\cdot\vect{u} \right)$.

To solve these equations, we use the OSMP numerical scheme, which
is a one-step high-order scheme developed in
\cite{Daru2004,Daru2009}. It is derived using a coupled space-time
Lax-Wendroff approach, where the formal order of accuracy in the
scalar case can be set at arbitrary order (in this paper, we use
order 11, that is the OSMP11 scheme). Imposing the MP conditions
(Monotonicity Preserving) prevents the appearance of numerical
oscillations in the vicinity of discontinuities while simultaneously
avoiding the numerical diffusion of extrema. In one space dimension,
on a uniform mesh with step-size $\Delta x$, at order $p$, it can be
written:
\begin{equation*}
  w_j^{n+1} = w_j^n-\frac{\Delta t}{\Delta x}(f_{j+\frac{1}{2}}^p-f_{j-\frac{1}{2}}^p)
\end{equation*}
where $f_{j+\frac{1}{2}}^p$ is the $p$th-order-accurate numerical
flux of the scheme at the cell interface $(j+\frac{1}{2})$.
Given the $l$
eigenvectors of the Jacobian matrix of the flux $\vect{r}_k$ and eigenvalues $\lambda_{k}$, the general
expression of the numerical fluxes can be written~:
\begin{equation}
  f_{j+\frac{1}{2}}^p = f_{j+\frac{1}{2}}^{\text{Roe}}+\frac{1}{2}\sum_k{(\psi^p\vect{r})_{k,j+\frac{1}{2}}}
\end{equation}
where, for clarity, the superscript $n$ has been
omitted. $\displaystyle f^{\text{Roe}}_{j+\frac{1}{2}}$ is
the first order Roe flux defined as follows:
\begin{equation}
f^{\text{Roe}}_{j+\frac{1}{2}} = \frac{1}{2} (f_{j} + f_{j+1}) - \frac{1}{2} \sum_k
( \delta  \vert f \vert \vect{r})_{k,j+\frac{1}{2}}
\end{equation}
with
\begin{equation*}
\delta \vert f \vert_{k,j+1/2} = \vert {\lambda} \vert_{k,j+1/2}
\delta {\alpha}_{k,j+1/2}
\end{equation*}
$\displaystyle \delta
{\alpha}_{k,j+\frac{1}{2}} = \vect{r}_{k} \cdot \left( w_{j+1}^n - w_{j}^n
\right) $ being the $k$-th Riemann invariant of the Jacobian matrix.
The $\psi^p$ are
corrective terms to obtain order $p$.
The function $\psi$ can be decomposed in odd and even parts:
\begin{equation} \label{OS_psi}
\psi^p_{k,j+\frac{1}{2}} = \sum_{n=1}^{ m }{\psi^{2 n}_{k,j+\frac{1}{2}}}+ js
\sum_{n=1}^{ m1 }{\psi^{2 n +1}_{k,j+1-\frac{js}{2}}} 
\end{equation} 
where $m=\lfloor \frac{p}{2} \rfloor$, $m1=\lfloor \frac{(p-1)}{2} \rfloor$
($\lfloor\;\rfloor$ is the integer division symbol), and
$js=\text{sign}(\lambda_{k,j+\frac{1}{2}})$. The odd and even $\psi$ functions
are given by the recurrence formulae (valid for $n \geq 1$):
\begin{equation} \label{OS_psi_even}
\psi^{2n}_{k,j+\frac{1}{2}} =\sum_{l=0}^{2n-2}{(-1)^l \mathbf{C}_{2n-2}^l
\cdot  (c_{2n} \delta {\alpha})_{k,j+\frac{1}{2}+n-1-l}} 
\end{equation}
\begin{equation} \label{OS_psi_odd} \psi^{2n+1}_{k,j+\frac{1}{2}}
=\sum_{l=0}^{2n-1}{(-1)^l \mathbf{C}_{2n-1}^l \cdot (c_{2n+1}  \delta
{\alpha})_{k,j+\frac{1}{2}+(n-1-l)\cdot js}}, 
\end{equation} 
where $\mathbf{C}_r^s=\frac{r!}{[(r-s)!s!]}$. The coefficients $c_q$
depend on the local CFL number, $\displaystyle {\nu}_{k,j+\frac{1}{2}}=
\frac{\delta t}{\delta x} {\lambda}_{k,j+\frac{1}{2}}$, and are given
by: 
\begin{equation} \label{cp}
  \begin{array}{l} 
    \displaystyle  (c_2)_{k,j+\frac{1}{2}}=\vert {\lambda}\vert_{k,j+\frac{1}{2}} (1-\vert \nu \vert_{k,j+\frac{1}{2}}) \\
    \displaystyle  (c_{q+1})_{k,j+\frac{1}{2}}=\frac{\vert\nu\vert_{k,j+\frac{1}{2}}+(-1)^{q}\lfloor \frac{(q+1)}{2} \rfloor}{q+1}
\cdot(c_q)_{k,j+\frac{1}{2}}, \ \ q \geq 2  \end{array}.
\end{equation}
At order $p$, the stencil of the scheme uses $p+2$ points. Flux limiting
TVD or MP constraints are then be applied to $\psi^p$ to make the scheme
non-oscillatory. The detail of the limiting procedure can be found in
\cite{Daru2009}. 

Near cut-cells, the existence of an adequate stencil of fluid points is
not necessarily provided for. Two main types of solutions can be
devised: either lower the order of accuracy and thus the stencil width,
or construct fictitious fluid values in the solid. We resort to the
second solution, with simple mirroring conditions with respect to the
solid boundary. The solution is satisfactory as long as the solid is
larger than the stencil of the scheme, which is the case for the
numerical examples considered in this paper. In case this condition
fails, we could resort to Ghost Fluid-type methods as in \cite{IBM25}.  

In two dimensions, the fluxes are computed using a directional
Strang splitting \cite{Strang} which is second-order accurate. However
the error of the scheme remains very low \cite{Daru2004}. This splitting
procedure will be expressed in section \ref{sec:description} devoted to
the coupling algorithm. 

\section{Treatment of the cells cut by the solid boundary in the
  Embedded Boundary method}
\label{sec:cutcells}

In this section, we recall the main ideas of the Embedded Boundary
method as exposed in \cite{IBM18,IBM30,IBM10}.

In order to take into account the position of the solid in the fluid
domain, we rely on the Embedded Boundary method, which consists in
modifying the fluid fluxes in cells that are cut by the solid
boundary (named cut cells), as in \cite{IBM10,IBM30}. At time $t$,
for a cut cell $\mathcal{C}$, we assume that the solid occupies a
volume fraction $\alpha_{\mathcal{C}}$. We also assume that the
density, velocity and pressure are constant in the cell. The fluid
mass, momentum and energy quantities contained in the cell are
therefore equal to their value at the center of the cell times the
volume of the cell and the volume fraction of fluid
$1-\alpha_{\mathcal{C}}$. In the same way, the computed fluxes are
assumed to be constant on the faces of a cell. Denoting by
$\kappa_{\mathcal{C}_1\mathcal{C}_2}$ the solid surface fraction of
the face between cells $\mathcal{C}_1$ and $\mathcal{C}_2$, the
effective flux between $\mathcal{C}_1$ and $\mathcal{C}_2$ is  the
computed flux times the surface of their interface times the fluid
surface fraction $1-\kappa_{\mathcal{C}_1\mathcal{C}_2}$. Additional
fluxes come from the presence of the moving solid boundary. These fluxes
arise due to the change in surface fractions and the work of the fluid
pressure on the solid surface. They are 
expressed in order to yield exact conservation of fluid mass and of
the total momentum and energy of the system.

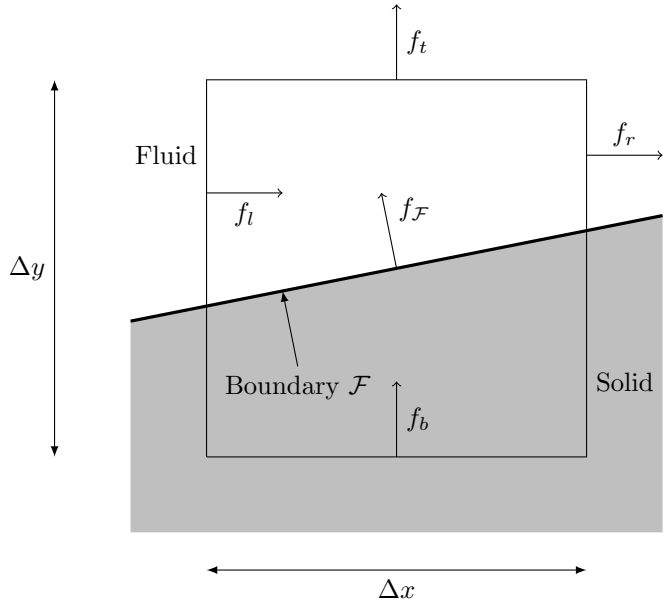
\begin{figure}[!ht]
  \centering
  \begin{tikzpicture}
    \filldraw[draw=none,fill=lightgray] (-1,-1) -- (6,-1) -- (6,3.2) -- (-1,1.8) -- (-1,-1);
    \draw (0,0) -- (0,5) -- (5,5) -- (5,0) -- (0,0);
    \draw[very thick] (-1,1.8) -- (6,3.2);
    \draw (5,1) node[anchor=west]{Solid};
    \draw (0,4) node[anchor=east]{Fluid};
    \draw[->] (0,3.5) -- (1,3.5);
    \draw (0.5,3.5) node[anchor=north]{$f_l$};
    \draw[->] (5,4) -- (6,4);
    \draw (5.5,4) node[anchor=south]{$f_r$};
    \draw[->] (2.5,5) -- (2.5,6);
    \draw (2.5,5.5) node[anchor=west]{$f_t$};
    \draw[->] (2.5,0) -- (2.5,1);
    \draw (2.5,0.5) node[anchor=west]{$f_b$};
    \draw[->] (2.5,2.5) -- (2.3,3.5);
    \draw (2.4,3) node[anchor=south west]{$f_{\mathcal{F}}$};
    \draw[latex-latex] (0,-1.5) -- (5,-1.5);
    \draw (2.5,-1.5) node[anchor=north]{$\Delta x$};
    \draw[latex-latex] (-2,0) -- (-2,5);
    \draw (-2,2.5) node[anchor=east]{$\Delta y$};
    \draw[-latex] (1.2,1.2) -- (1,2.2);
    \draw (1.2,1.2) node[anchor=north]{Boundary $\mathcal{F}$};
  \end{tikzpicture}
  \caption{Physical fluxes in a cut cell}\label{fig:2Dphys}
\end{figure}

For the sake of simplicity, we limit ourselves to two space
dimensions. However the three-dimensional case can be carried out in
a similar way. Let us consider a fluid cell $\mathcal{C}$ cut by the
boundary, as shown in Figure \ref{fig:2Dphys}. The indices $l$, $r$,
$t$ and $b$ indicate respectively left, right, top and bottom in the
sequel.
Integrating the Euler equations on the cut cell and over the time
interval $[ n \Delta t, (n+1) \Delta t ]$, and applying the divergence
theorem, we get:
  \begin{multline}
    (1-\alpha_{\mathcal{C}}^{n+1})\Delta w_{\mathcal{C}} = \Delta
    t\left(\frac{1-\kappa_{\mathcal{C}l}}{\Delta
    x}{f}_{\mathcal{C}l} -
    \frac{1-\kappa_{\mathcal{C}r}}{\Delta
    x}{f}_{\mathcal{C}r} +
    \frac{1-\kappa_{\mathcal{C}b}}{\Delta
    y}{f}_{\mathcal{C}b} \right.\\
    \left.-
    \frac{1-\kappa_{\mathcal{C}t}}{\Delta
    y}{f}_{\mathcal{C}t} \right)  +
    \frac{\Delta t}{\Delta x\Delta y}X_{\mathcal{F}} +
    \sum_{\mathcal{F}\in\mathcal{C}}{\Delta w_{\mathcal{F}}^n}
    \label{eqn:fluxmodif_new}
  \end{multline}
where $\Delta
w_{\mathcal{C}}=w_{\mathcal{C}}^{n+1}-w_{\mathcal{C}}^{n}$ is the
time increment and all fluxes are time-averaged over the time
interval (the time averaging will be specified later). At the solid
walls, pressure forces cause momentum and energy exchange between
the solid and the fluid. They are taken into account through the
exchange term $X_{\mathcal{F}}$. The detailed expression of
$X_{\mathcal{F}}$ will be given in section \ref{sec:step5}. Finally,
the quantity $\Delta w_{\mathcal{F}}^n$ represents the amount of $w^n$
swept by each solid boundary $\mathcal{F}$ present in the cell during the time
step. The solid boundary $\mathcal{F}$ is the largest subsegment of the solid boundary which is contained in one single cell (not necessarily the same) at times $n\Delta t$ and $(n+1)\Delta t$. The precise definition of $\mathcal{F}$ and the expression of $\Delta w_{\mathcal{F}}^n$ will be given in section \ref{sec:step4}.

\section{Coupling algorithm}
\label{sec:description}

Since the Discrete Element method is computationally
expensive, the coupling algorithm should be explicit in order to avoid
costly iterative procedures. In fact, the CFL condition of the explicit
time-scheme gives the appropriate criterion for the capture of the
high-frequency eigenmodes involved in the solid body fast
dynamics. Moreover, as it is well known, explicit methods are more
robust for impact problems.  We choose the following general structure of the
algorithm, which can be traced back to Noh \cite{Noh}:
\begin{itemize}
\item The position of the solid and the density, velocity and pressure
  of the fluid are known at time $t$
\item The fluid exerts a pressure force on the solid boundaries:
  knowing the total forces applied on the solid, the position of the
  solid is advanced to time $t+\Delta t$
\item The density, velocity and pressure of the fluid are then
  computed at time $t+\Delta t$. This step takes into account the new
  position and velocity of the solid boundary, as well as the work of
  the forces of pressure on the boundary during the time step.
\end{itemize}

The choice of the coupling algorithm is guided by the conservation of
the global momentum and energy of the system and the conservation of
constant flows (see section \ref{sec:constant}).

At the beginning of a time step, at time $n\Delta t$, the position and
rotation of the solid particle $(\vect{X}^n,\matr{Q}^n)$, the velocity
and angular velocity of the solid particle
$(\vect{V}^n,\vect{\Omega}^n)$ and the fluid state $\vect{w}^n$ are known.
We choose the following general architecture for the algorithm:

\begin{center}
  \resizebox{1.1\textwidth}{!}{
    \begin{tikzpicture}
      \centering
      \draw (0,0.5) node[anchor=north]{SOLID};
      \draw[rounded corners] (-1,0.5) -- (1,0.5) -- (1,0) -- (-1,0) --
      cycle;
      \draw (7,0.5) node[anchor=north]{FLUID};
      \draw[rounded corners] (6,0.5) -- (8,0.5) -- (8,0) -- (6,0) --
      cycle;
      \draw (3.5,0.5) node[anchor=north]{COUPLING};
      \draw[rounded corners] (2.25,0.5) -- (4.75,0.5) -- (4.75,0) -- (2.25,0) --
      cycle;
      \draw (0,-1) node[anchor=north]{$\vect{X}^n$, $\matr{Q}^n$,
    $\vect{V}^n$, $\matr{P}^n$};
      \draw (6.75,-1) node[anchor=north]{$\rho^n$, $\vect{u}^n$, $p^n$};
      \draw[rounded corners] (-2,-1) -- (2,-1) -- (2,-1.7) -- (-2,-1.7) --
      cycle;
      \draw[rounded corners] (5.5,-1) -- (8,-1) -- (8,-1.7) -- (5.5,-1.7) --
      cycle;
      \draw[-latex,dashed] (6,-1.7) -- (6,-4.4);
      \draw (6,-2.7) node[anchor=east]{
        \begin{minipage}{2.75cm}
          (1) Computation of fluxes $f$
        \end{minipage}};
      \draw (6,-4.4)  node[anchor=north]{$\overline{p}_x$,
        $\overline{p}_y$, $f$};
      \draw[rounded corners] (5,-4.4) -- (7,-4.4) -- (7,-5.1) -- (5,-5.1) -- cycle;
      \draw[-latex] (5,-4.75) -- (0,-4.75);
      \draw (2.5,-4.75) node[anchor=north]{
    \begin{minipage}{4cm}
          (2) Predicted pressure is transferred to the solid boundary
      \end{minipage}};
      \draw[-latex] (0,-1.7) -- (0,-7.1);
      \draw (0,-4.4) node[anchor=east]{
    \begin{minipage}{3cm}
          (3) Solid step (using predicted boundary pressure)
      \end{minipage}};
      \draw (-0.5,-7.1) node[anchor=north]{$\vect{X}^{n+1}$,
    $\matr{Q}^{n+1}$, $\vect{V}^{n+1}$, $\matr{P}^{n+1}$};
      \draw[rounded corners] (2,-7.1) -- (2,-8.) -- (-3,-8.) -- (-3,-7.1) -- cycle;
      \draw[-latex] (2,-7.55) -- (7.5,-7.55);
      \draw (4.75,-7.55) node[anchor=north]{
        \begin{minipage}{4.5cm}
          (4) Update of the boundary position, computation of the $\alpha^{n+1}$ and $\kappa^{n+1}$
      \end{minipage}};
      \draw[-latex] (7.5,-1.7) -- (7.5,-8.7);
      \draw (7.5,-5.2) node[anchor=west]{
    \begin{minipage}{6cm}
          (5) Fluid update:\\
          \begin{equation*}
            \left\{
            \begin{array}{l}
              \rho^{n+1}=\rho^n+\Delta\rho \\
              \rho^{n+1}\vect{u}^{n+1}=\rho^n\vect{u}^n+\Delta
              (\rho\vect{u}) \\
              \rho^{n+1}E^{n+1}=\rho^n E^n+\Delta (\rho E)
            \end{array}
            \right.
          \end{equation*}
      \end{minipage}};
      \draw (7.5,-8.7) node[anchor=north]{$\rho^{n+1}$, $\vect{u}^{n+1}$, $p^{n+1}$};
      \draw[rounded corners] (6,-8.7) -- (6,-9.6) -- (9,-9.6) -- (9,-8.7) -- cycle;
    \end{tikzpicture}
  }
\end{center}

Steps (1) to (5) of the algorithm are computed successively, and are
detailed in the following subsections.

\subsection{Computation of fluid fluxes and of the boundary pressure (steps (1) and (2))}

Step (1) is a precomputation of fluxes without considering the
presence of a solid boundary. As said above, the fluxes are computed
in every cell using the OSMP11 scheme. However, we emphasize that
the coupling algorithm does not depend on the choice of the
numerical scheme. The fluxes are then stored for later use in step
(5).

The other aim of this step is the computation of mean pressures in
each cut-cell during the time-step in each direction
$\overline{p}_x$ and $\overline{p}_y$. These pressures, transferred
to the solid boundary in step (2),  account for the forces exerted
by the fluid on the solid during the time-step. The same mean
pressures will be used in step (5) to compute the momentum and
energy exchanged between the solid and the fluid. In this way, the
choice of $\overline{p}_x$ and $\overline{p}_y$ has no effect on the
conservation of fluid mass, momentum or energy of the system. On the
contrary it is a key ingredient for the exact conservation of
constant flows (see section \ref{sec:constant}). The explicit
structure of our solid and fluid methods allows several
possibilities for the choice of boundary pressures while maintaining
the stability of the coupling algorithm. This is unusual in
fluid-structure interaction.

The Strang directional splitting algorithm \cite{Strang} is originally formulated as follows:
\begin{equation*}
  \vect{w}_j^{(n+1)} = L_x\left(\frac{\Delta t}{2}\right) L_y\left(\Delta t\right) L_x\left(\frac{\Delta t}{2}\right)\vect{w}_j^n ,
\end{equation*}
where $L_x(\Delta t)$ and $L_y(\Delta t)$ are finite-difference
approximation operators for the integration by a time-step $\Delta t$ in
directions $x$ and $y$ respectively. 
Here, this splitting procedure is implemented in a simplified form: 
\begin{equation*}
	\vect{w}_j^{(n+2)} = L_x\left(\Delta t\right) L_y\left(\Delta t\right) L_y\left(\Delta t\right) L_x\left(\Delta t\right)\vect{w}_j^n ,
\end{equation*}
that recovers the symmetry of the solution every two time steps.
In our case, $L_x$ and $L_y$
involve the computation of a flux in the $x$ or $y$ direction using the
state of fluid $\vect{w}$ of the cells. The mean pressures
$\overline{p}_x$ and $\overline{p}_y$ are then the pressures in the cell
used for the computation of the fluxes by $L_x$ and $L_y$, respectively. An
analogous definition could be derived for other time integration
methods, such as Runge-Kutta for instance. The directional splitting
used for the fluid flux computation does not require the solid body
displacement to be split in $x$ and $y$ components. It is applied here
only to recover second-order accuracy of the fluxes.

\subsection{Computation of the solid step (step (3))}
 Step (3) consists mainly in the application of the time integration
 scheme for the rigid body motion described in section \ref{sec:solide}. The essential difference with
 an uncoupled version lies in the integration of boundary pressure forces.
 As we consider an explicit coupling, the only
 boundary pressures available are $\overline{p}_x$ and
 $\overline{p}_y$.

 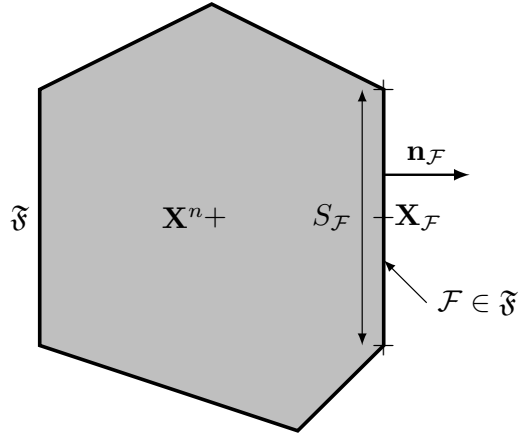
\begin{figure}[!ht]
   \centering
   \resizebox{0.45\textwidth}{!}{
     \begin{tikzpicture}
       \draw[fill=lightgray,very thick] (2cm,-2cm) -- (2cm,1cm) -- (0cm,2cm) -- (-2cm,1cm) -- (-2cm,-2cm) -- (1cm,-3cm) --
       cycle;
       \draw (0.0416667cm,-0.5cm) node{+};
       \draw (0.0416667cm,-0.5cm) node[anchor = east] {$\vect{X}^n$};
       \draw [latex-latex] (1.75cm,-2cm) -- (1.75cm,1cm) ;
       \draw (1.75cm,-0.5cm) node[anchor=east]{$S_{\mathcal{F}}$};
       \draw [-latex,thick] (2cm,0cm) -- (3cm,0cm);
       \draw (2.5cm,0) node[anchor=south]{$\vect{n}_{\mathcal{F}}$};
       \draw (2cm,-0.5cm) node{+};
       \draw (2cm,-0.5cm) node[anchor=west]{$\vect{X}_{\mathcal{F}}$};
       \draw (-2cm,-0.5cm) node[anchor=east]{$\mathfrak{F}$};
       \draw (2cm,-2cm) node{+};
       \draw (2cm,1cm) node{+};
       \draw [-latex] (2.5cm,-1.5cm) -- (2cm,-1cm);
       \draw (2.5cm,-1.5cm) node[anchor=west]{$\mathcal{F}\in\mathfrak{F}$};
     \end{tikzpicture}
   }
   \caption{Geometric description of the particles}\label{fig:particles}
 \end{figure}

 The solid is assumed to be polygonal (in two space dimensions) as
 described in Figure \ref{fig:particles}.
 We denote by $\mathfrak{F}$ the list of all faces of the solid in
 contact with fluid. For every face $\mathcal{F}\in\mathfrak{F}$, the
 position of
 the center of the face is given by vector $\vect{X}_{\mathcal{F}}$,
 and we denote by $S_{\mathcal{F}}$ its surface and
 $\vect{n}_{\mathcal{F}}$ its  normal vector
 (oriented from the solid to the fluid). The fluid pressure force
 $\vect{F}_{\mathcal{F}}$ exerted on face $\mathcal{F}\in\mathfrak{F}$ is
 written as:
 \begin{align}
   \vect{{F}}_{\mathcal{F}}\cdot\vect{e}_x &=
   -\overline{p}_xS_{\mathcal{F}}n^x_{\mathcal{F}} \label{eqn:forcepressionx}\\
   \vect{{F}}_{\mathcal{F}}\cdot\vect{e}_y &=
   -\overline{p}_yS_{\mathcal{F}}n^y_{\mathcal{F}}\label{eqn:forcepressiony}
 \end{align}
 The total fluid pressure force $\vect{{F}}_{f}^n$ is the sum of the
 contributions on each face:
 \begin{equation}
   \vect{F}_{f}^n =
   \sum_{\mathcal{F}\in\mathfrak{F}}{\vect{{F}}_{\mathcal{F}}}\label{eqn:forcepression}
 \end{equation}
 The fluid pressure torque $\vect{\mathcal{M}}_{f}^n$ is the sum of the torques
 of the pressure forces at the center of mass of the solid body:
 \begin{equation*}
   \vect{\mathcal{M}}_{f}^n =
   \sum_{\mathcal{F}\in\mathfrak{F}}{\vect{{F}}_{\mathcal{F}}\wedge\left(\vect{X}^n-\vect{X}_{\mathcal{F}}
     \right)}
 \end{equation*}

 The solid time-step is written as in equations (\ref{eqn:vitesse1})
 to (\ref{eqn:contrainte2}), with the only difference that the fluid
 pressure force and torque are taken constant during the whole
 time-step, equal to $\vect{F}_f^n$ and $\vect{\mathcal{M}}_f^n$
 (including in equations (\ref{eqn:vitesse2}) and
 (\ref{eqn:vitesserot2})). The fact that $\vect{F}_f^n$,
 $\vect{\mathcal{M}}_f^n$, $\vect{V}^{n+\frac{1}{2}}$ and
 $\matr{P}^{n+\frac{1}{2}}$ are constant during the time-step will be
 used in the conservation analysis in section \ref{sec:conservation}.

\subsection{Update of the boundary and of the volume fractions (step (4))} \label{sec:step4}

Several tasks are carried out in step (4). For each cell
$\mathcal{C}$, the new solid volume fraction of the cell
$\alpha_{\mathcal{C}}^{n+1}$ and new surface fractions
$\kappa_{\mathcal{C}}^{n+1}$ are computed. In addition, for each
solid boundary $\mathcal{F}$, the pressures $\overline{p}_x$ and
$\overline{p}_y$ are stored and the swept quantities
$\Delta w_{\mathcal{F}}^n$ used in
(\ref{eqn:fluxmodif_new})   are evaluated.

\begin{figure}[!ht]
  \centering
  \resizebox{0.8\textwidth}{!}{
    \begin{tikzpicture}[x=2cm,y=2cm]
      \draw[fill=lightgray,draw=none] (0.613692,0.886308) -- (1,1.375)
      -- (1.437056,1.0272) -- (1,0.5) -- cycle;
      \draw (-0.5,0) -- (3.5,0);
      \draw (-0.5,1) -- (3.5,1);
      \draw (-0.5,2) -- (3.5,2);
      \draw (0,-0.5) -- (0,2.5);
      \draw (1,-0.5) -- (1,2.5);
      \draw (2,-0.5) -- (2,2.5);
      \draw (3,-0.5) -- (3,2.5);
      \draw[thick] (-0.6,2.1) -- (2.1,-0.6);
      \draw[thick] (-0.5,2.5) -- (3.5,-0.5);
      \draw (-0.55,1.95) -- (-0.45,2.05);
      \draw (-0.05,1.45) -- (0.05,1.55);
      \draw (0.45,0.95) -- (0.55,1.05);
      \draw (0.95,0.45) -- (1.05,0.55);
      \draw (1.45,-0.05) -- (1.55,0.05);
      \draw (1.95,-0.55) -- (2.05,-0.45);
      \draw (-0.03,2.085) -- (0.03,2.165);
      \draw (0.13666,1.96) -- (0.19666,2.04);
      \draw (0.97,1.335) -- (1.03,1.415);
      \draw (1.47,0.96) -- (1.53,1.04);
      \draw (1.97,0.585) -- (2.03,0.665);
      \draw (2.80333,-0.04) -- (2.86333,0.04);
      \draw (2.97,-0.165) -- (3.03,-0.085);
      \draw (-0.29,2.26) -- (-0.23,2.34);
      \draw[dotted] (-0.5,2) -- (-0.26,2.3);
      \draw (0.275685,1.835735) -- (0.335685,1.915735);
      \draw[dotted] (0,1.5) -- (0.305685,1.875735);
      \draw (0.84137,1.41147) -- (0.90137,1.49147);
      \draw[dotted] (0.5,1) -- (0.87137,1.45147);
      \draw (1.407056,0.9872) -- (1.467056,1.0672);
      \draw[dotted] (1,0.5) -- (1.437056,1.0272);
      \draw (1.97274,0.5629437) -- (2.03274,0.6429437);
      \draw[dotted] (1.5,0) -- (2.00274,0.6029437);
      \draw (2.5384,0.13868) -- (2.5984,0.21868);
      \draw[dotted] (2,-0.5) -- (2.5684,0.17868);
      \draw (-0.32019,1.721019) -- (-0.22019,1.82019);
      \draw[dotted] (-0.27019,1.77019) -- (0,2.125);
      \draw (-0.172877,1.572877) -- (-0.072877,1.672877);
      \draw[dotted] (-0.122877,1.622877) -- (0.1666,2);
      \draw (0.563692,0.836308) -- (0.663692,0.936308);
      \draw[dotted] (0.613692,0.886308) -- (1,1.375);
      \draw (1.005633,0.394367) -- (1.105633,0.494367);
      \draw[dotted] (1.055633,0.444367) -- (1.5,1);
      \draw (1.44757,-0.0475745) -- (1.54757,0.0524255);
      \draw[dotted] (1.49757,0.0024255) -- (2,0.625);
      \draw (0.806846,0.693154) node[anchor= north
      east,fill=white,rounded corners=3pt,draw=black]{$\Phi_n(\mathcal{F})$};
      \draw (1.218528,1.2011) node[anchor= south
      west,fill=white,rounded corners=3pt,draw=black]{$\Phi_{n+1}(\mathcal{F})$};
      \draw[very thick] (0.613692,0.886308) -- (1,0.5);
      \draw[very thick] (1,1.375) -- (1.437056,1.0272);
      \draw (1.75,-0.25) node[anchor=north east,fill=white,rounded corners=3pt,draw=black]{Solid boundary at
        $t=n\Delta t$};
      \draw (0.5,1.75) node[anchor=south west,fill=white,rounded corners=3pt,draw=black]{Solid boundary at
        $t=(n+1)\Delta t$};
      \draw[-latex] (2.2,0.8) -- (1.1,0.9);
      \draw (2.2,0.8) node[anchor=west,fill=white,rounded corners=3pt,draw=black]{Domain of integration of $\Delta
        w^n_{\mathcal{F}}$};
    \end{tikzpicture}
  }
  \caption{Update of the boundary and computation of the $\Delta w^n_{\mathcal{F}}$}
  \label{fig:decoupage}
\end{figure}
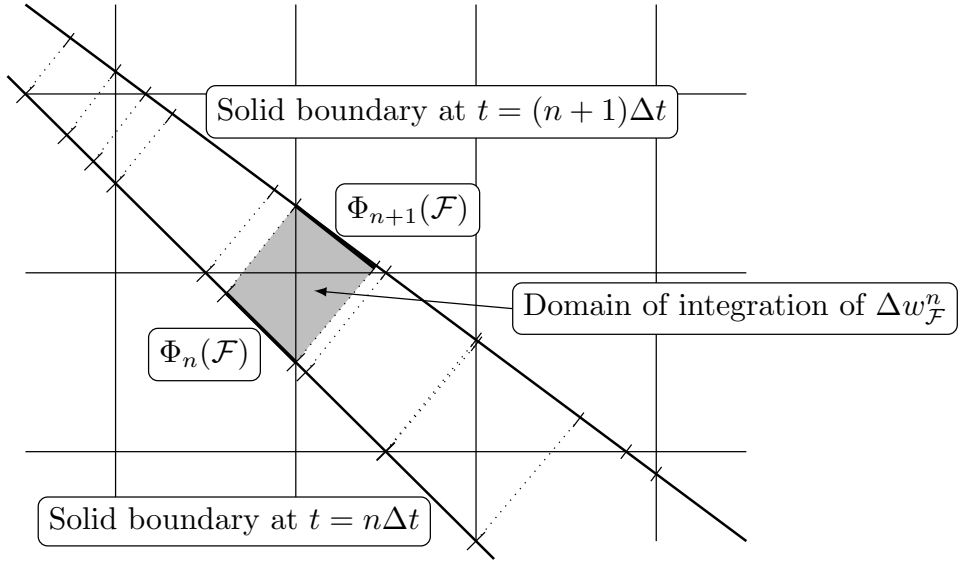

 In two dimensions, the solid boundary is
polygonal, and we therefore only have to deal with plane boundaries
$\mathcal{F}$. In order to simplify the computation of the average
of $\overline{p}_x$ and $\overline{p}_y$ on $\mathcal{F}$, we also
assume that each boundary $\mathcal{F}$ is contained only in one
cell at time $n\Delta t$. The computation of the contribution of
$\Delta w_{\mathcal{F}}^n$ to each cell is also easier if
$\mathcal{F}$ is entirely in the cell at time $(n+1)\Delta t$. We
denote $\Phi_n(\mathcal{F})$ the position of boundary $\mathcal{F}$
at time $n\Delta t$. We choose to define $\mathcal{F}$ as the largest
subsegment of the boundary polygon such that $\Phi_n(\mathcal{F})$ is
contained in cell $\mathcal{C}_n$ at time $n\Delta t$ and
$\Phi_{n+1}(\mathcal{F})$ is contained in cell $\mathcal{C}_{n+1}$ at
time $(n+1)\Delta t$ (see Figure \ref{fig:decoupage}). The two cells
need not be necessarily different. $\mathcal{F}$ may contain one or both
vertices of the polygonal boundary at its ends, but we assume that each
$\mathcal{F}$ is contained in one single polygonal face. At each new
time step, the polygonal boundary is subdivised into a new set of plane
boundaries $\mathcal{F}$. 
Each newly computed boundary $\mathcal{F}\in\mathfrak{F}$ stores
every variable necessary for the coupling: the surface
$S_{\mathcal{F}}$ and the normal vector $\vect{n}_{\mathcal{F}}$ of
$\Phi_n(\mathcal{F})$, the center of mass $\vect{X}_{\mathcal{F}}$
of $\mathcal{F}$, and we define
$\vect{X}_{\mathcal{F}}^0=\Phi_0(\vect{X}_{\mathcal{F}})$. The
boundary also stores the pressures $\overline{p}_x$ and
$\overline{p}_y$ in the cell occupied by $\Phi_n(\mathcal{F})$, and
the velocity of the center of the boundary at time
$(n+\frac{1}{2})\Delta t$, $\vect{V}_{\mathcal{F}}^{n+\frac{1}{2}}$,
computed as:
  \begin{equation}
    \vect{V}_{\mathcal{F}}^{n+\frac{1}{2}} =
    \vect{V}^{n+\frac{1}{2}} +
    \matr{P}^{n+\frac{1}{2}}\matr{D}^{-1}\cdot\left(\vect{X}^0-\vect{X}_{\mathcal{F}}^0\right) \label{eqn:vitesseface}
  \end{equation}

The swept quantities $\Delta w_{\mathcal{F}}^n$ are computed as the
integral of $w^n$ in the quadrangle bounded by $\Phi_n(\mathcal{F})$
and $\Phi_{n+1}(\mathcal{F})$ (see Fig. \ref{fig:decoupage}). The condition
\begin{equation}
  \sum_{\mathcal{F}}{\Delta w_{\mathcal{F}}^n} =
  \sum_{\mathcal{C}}{(\alpha_{\mathcal{C}}^{n+1}-\alpha_{\mathcal{C}}^{n})w_{\mathcal{C}}^n}
  \label{eqn:balayage}
\end{equation}
is then automatically satisfied
as the set of such quadrangles is a partition of the volume swept by
the solid during the time step.

The computation of $\alpha_{\mathcal{C}}^{n+1}$ and
$\kappa_{\mathcal{C}}^{n+1}$ involves the intersection of planes with
rectangles, and can be carried out geometrically. The application of the
divergence theorem to cell
$\mathcal{C}$ shows that the following relations are satisfied:
\begin{align}
  \kappa_{\mathcal{C}l}^{n} &= \kappa_{\mathcal{C}r}^{n} +
  \sum_{\mathcal{F}\in\mathcal{C}}{\frac{S_{\mathcal{F}}}{\Delta
  y}n_{\mathcal{F}}^x} \label{eqn:compatibilitex} \\
  \kappa_{\mathcal{C}b}^{n} &=
  \kappa_{\mathcal{C}t}^{n} +
  \sum_{\mathcal{F}\in\mathcal{C}}{\frac{S_{\mathcal{F}}}{\Delta
  x}n_{\mathcal{F}}^y} \label{eqn:compatibilitey}
\end{align}
These conditions correspond to the geometric conservation laws (GCL)
in ALE methods \cite{LesoinneFarhat}, and will be used in the
analysis of the consistency of our method.

\subsection{Modification of fluxes (step (5))}\label{sec:step5}

Step (5) of the algorithm consists mainly in the  computation of the
final values $w^{n+1}_{\mathcal{C}}$ in each cell, using a fully
discrete expression of Eq. (\ref{eqn:fluxmodif_new}). It is the only
part of the algorithm where the fluid ``sees'' the presence of a
solid. The explicit fluid fluxes were pre-computed in step (1) on
the Cartesian regular grid, and the modification of fluxes aims at
conserving the mass of fluid and balancing the momentum and energy
transferred to the solid during the time-step.

The exchange term in (\ref{eqn:fluxmodif_new}) can be written as
$$X_{\mathcal{F}}=\sum_{\mathcal{F}\in\mathcal{C}}{S_{\mathcal{F}}{f}_{\mathcal{F}}},$$
where ${f}_{\mathcal{F}}$ is the fluid flux at the solid boundary
$\mathcal{F}$ (see fig. \ref{fig:2Dphys}), that are approximated as:

  \begin{equation}
    {f}_{\mathcal{F}} =
      \frac{1}{S_{\mathcal{F}}}\left( 0,
      \int_{\mathcal{F}}{\overline{p}_xn_{\mathcal{F}}^x},
      \int_{\mathcal{F}}{\overline{p}_yn_{\mathcal{F}}^y},
     \vect{V}_{\mathcal{F}}^{n+\frac{1}{2}}\cdot\int_{\mathcal{F}}{\left(
        \begin{array}{c}
          \overline{p}_xn_{\mathcal{F}}^x \\ \overline{p}_yn_{\mathcal{F}}^y
        \end{array}
        \right) }
        \right)^T \label{eqn:fluxparoisolide}
  \end{equation}
 Here $\vect{V}_{\mathcal{F}}^{n+\frac{1}{2}}$ is the velocity of the
  center of the boundary and is defined in (\ref{eqn:vitesseface}).

Using the fluid fluxes given by the OSMP11 scheme, we finally
compute the time increment $\Delta w_{\mathcal{C}}$ from the
following fully discrete version of equation
(\ref{eqn:fluxmodif_new}):
  \begin{multline}
    (1-\alpha_{\mathcal{C}}^{n+1})\Delta w_{\mathcal{C}} = \Delta
    t\left(\frac{1-\kappa_{\mathcal{C}l}^{n+1}}{\Delta
    x}{f}_{\mathcal{C}l} -
    \frac{1-\kappa_{\mathcal{C}r}^{n+1}}{\Delta
    x}{f}_{\mathcal{C}r} +
    \frac{1-\kappa_{\mathcal{C}b}^{n+1}}{\Delta
    y}{f}_{\mathcal{C}b} \right.\\
    \left.-
    \frac{1-\kappa_{\mathcal{C}t}^{n+1}}{\Delta
    y}{f}_{\mathcal{C}t} \right)  +
    \frac{\Delta t}{\Delta x\Delta y}\sum_{\mathcal{F}\in\mathcal{C}}{S_{\mathcal{F}}{f}_{\mathcal{F}}} +
    \sum_{\mathcal{F}\in\mathcal{C}}{\Delta w^n_{\mathcal{F}}}
    \label{eqn:fluxmodif2}
  \end{multline}
The value of $w_{\mathcal{C}}$ is then updated in every cell:
   $\displaystyle w^{n+1}_{\mathcal{C}} = w^{n}_{\mathcal{C}} + \Delta
    w_{\mathcal{C}}$. 

A main difference with \cite{IBM30} lies in the time integration of
cell-face apertures $(1-\kappa_{\mathcal{C}})$. Falcovitz {\it et
al.} \cite{IBM30} use time-averaged cell-face apertures over the
time step (at time $(n+\frac{1}{2})\Delta t$), ensuring consistency
(in the sense that the uniform motion of a solid-fluid system is
exactly preserved). In fact, a key ingredient in the consistency proof is the fact that conditions (\ref{eqn:compatibilitex}) and (\ref{eqn:compatibilitey}) are checked exactly. In \cite{IBM30}, the consistent choice of time-averaged cell-face apertures $\kappa^{n+\frac{1}{2}}$ and solid surfaces in the fluid cell $\tilde{S}^{n+\frac{1}{2}}$ allows to check these conditions.

Here we instead take $\kappa^{n+1}$ and recover
consistency using the solid surface $\tilde{S}^{n+1}$ present in the
fluid cell at time $(n+1)\Delta t$. This result is proved in section
\ref{sec:constant}. Note that
$\tilde{S}_{\mathcal{C}}^{n+1}=\sum_{\Phi_{n+1}(\mathcal{F})\in\mathcal{C}}{S_{\mathcal{F}}^n}$
as the solid is undeformable. This choice is motivated by the fact
that the computation of time averaged $\kappa^{n+\frac{1}{2}}$ and
$\tilde{S}^{n+\frac{1}{2}}$ is already complex in two dimensions and
might become intractable in three dimensions. In addition, it requires
an implicit resolution of $\tilde{S}^{n+\frac{1}{2}}$ in order to
preserve the energy of the system. The choice of $\kappa^{n+1}$
theoretically reduces the accuracy of the method in
cut-cells. However, the accuracy of our numerical results did not
advocate the use of the time-averaged $\kappa^{n+\frac{1}{2}}$ and the
related added complexity in the algorithm.

In order to avoid the classical restriction of the time-step due to
vanishing volumes:
\begin{equation*}
  \Delta t\leq\frac{(1-\alpha_{\mathcal{C}})\text{min}(\Delta x,\Delta
  y)}{\lVert\vect{u}\rVert+c},
\end{equation*}
where $c$ is the local speed of sound, we resort to the mixing of
small cut cells with their neighbors to prevent instabilities.

\subsection{Conservative mixing of small cut cells}
Two main methods have been developed to ensure the stability of
conservative Embedded Boundary methods. A first method consists in
computing a reference state using nonconservative interpolations,
modified by redistributing the conservation error on neighbouring cells
\cite{IBM18,IBM9,IBM27}. A second method is to compute a fully conservative
state using a formula similar to (\ref{eqn:fluxmodif2}). For stability
reasons, small cells are merged with neighboring cells
using a conservative procedure (originating from Glimm's idea
\cite{IBM28}). We choose this second class of method, as \cite{IBM10}
and \cite{IBM30}.

Target cells need to be defined for small cells to be merged with
them. \cite{IBM10} defines an equivalent normal vector to the boundary
in the cell and mixes the cells preferentially in that
direction. \cite{IBM30} rather merges newly exposed or newly covered
cells with full neighbours having a face in common. In order to deal with
cells occupied by several boundaries (impact of two solids), we cannot
define a normal vector in every cell and we choose to improve the
strategy applied in \cite{IBM30}. We define
small cells as $\alpha_{\mathcal{C}} > 0.5$. For mixing two cells
$\mathcal{C}$ and $\mathcal{C}_T$, so they have equal final value
$w$, the following quantities are exchanged:
\begin{align*}
  M_{\mathcal{C}\mathcal{C}_T} &=
  \frac{\alpha_{\mathcal{C}_T}}{\alpha_{\mathcal{C}}+\alpha_{\mathcal{C}_T}}\left(w_{\mathcal{C}_T}-w_{\mathcal{C}}\right)
  \\ M_{\mathcal{C}_T\mathcal{C}} &=
  \frac{\alpha_{\mathcal{C}}}{\alpha_{\mathcal{C}}+\alpha_{\mathcal{C}_T}}\left(w_{\mathcal{C}}-w_{\mathcal{C}_T}\right)
\end{align*}
and it is easy to check that
$w_{\mathcal{C}}+M_{\mathcal{C}\mathcal{C}_T}=w_{\mathcal{C}_T}+M_{\mathcal{C}_T\mathcal{C}}$.
In the two dimensional case, we select the target cell
$\mathcal{C}_T$ as the fully-fluid cell ($\alpha_{\mathcal{C}_T}=0$)
nearest to $\mathcal{C}$, such that the path between the two cells
does not cross a solid boundary. A recursive subroutine finds such a
target cell in a small number of iterations, without any restriction
on the geometry of the fluid domain.


\section{Analysis of the conservation of mass, momentum and energy}
\label{sec:conservation}

In this section, we analyze the conservation properties of the
coupling algorithm. These properties are verified for periodic boundary
conditions or for an infinite domain.

\subsection{Integration on the fluid domain}
\label{sec:fluide}

Integrating $w$ on the fluid domain $\Omega_f^{n+1}$ at time
$(n+1)\Delta t$, we obtain using (\ref{eqn:fluxmodif2}) and the
cancellation of fluxes on each cell face:
\begin{align*}
    \frac{1}{\Delta x\Delta y}\int_{\Omega_f^{n+1}}{w^{n+1}} =&
    \sum_{\mathcal{C}}{(1-\alpha_{\mathcal{C}}^{n+1})w_{\mathcal{C}}^{n}}
    + \sum_{\mathcal{C}}{(1-\alpha_{\mathcal{C}}^{n+1})\Delta
    w_{\mathcal{C}}^n}\\
  =& \sum_{\mathcal{C}}{(1-\alpha_{\mathcal{C}}^{n+1})w_{\mathcal{C}}^{n}}
    +
    \sum_{\mathcal{F}}{\frac{\Delta tS_{\mathcal{F}}}{\Delta
    x\Delta y}{f}_{\mathcal{F}}} + \sum_{\mathcal{F}}{\Delta
    w^n_{\mathcal{F}}}
\end{align*}

Using (\ref{eqn:balayage})  we finally get:
\begin{align}
  \frac{1}{\Delta x\Delta y}\int_{\Omega_f^{n+1}}{w^{n+1}} &=
  \sum_{\mathcal{C}}{(1-\alpha_{\mathcal{C}}^{n})w_{\mathcal{C}}^{n}}
  + \sum_{\mathcal{F}}{\frac{\Delta tS_{\mathcal{F}}}{\Delta x\Delta
  y}{f}_{\mathcal{F}}} \nonumber
\\ & = \frac{1}{\Delta x\Delta
  y}\int_{\Omega_f^{n}}{w^{n}} +  \frac{\Delta t}{\Delta x\Delta
  y}\sum_{\mathcal{F}}{S_{\mathcal{F}}{f}_{\mathcal{F}}} \label{eqn:conservfluid}
\end{align}
the expression of ${f}_{\mathcal{F}}$ being given in Eq.
\ref{eqn:fluxparoisolide}.

The first component of system (\ref{eqn:conservfluid}) expresses the
fluid mass conservation. In order to proceed with the analysis of
momentum and energy conservation, let us now turn to the solid part.

\subsection{Solid conservation balance}

Since the solid is treated using a Lagrangian method, the conservation
of solid mass is straightforward. The fluid pressure force applied
on the solid during the time step is given by
(\ref{eqn:forcepressionx}), (\ref{eqn:forcepressiony}) and
(\ref{eqn:forcepression}). Let us consider a solid boundary
$\mathcal{F}\in\mathfrak{F}$, and denote by $\Delta
\vect{\mathcal{P}}_{\mathcal{F}}$ the solid momentum variation
induced by the pressure forces on $\mathcal{F}$, and $\Delta
\mathcal{E}_{\mathcal{F}}$ the corresponding energy variation.
Recalling that the pressure forces are kept constant during the time
step, the balance of momentum and energy is given by:
\begin{align*}
  \Delta \vect{\mathcal{P}}_{\mathcal{F}} &= \Delta t
  \vect{{F}}_{\mathcal{F}} \\
  \Delta \mathcal{E}_{\mathcal{F}}
  &= \Delta t \vect{{F}}_{\mathcal{F}}\cdot\left(\frac{1}{S_{\mathcal{F}}}\int_{\mathcal{F}}{\vect{V}^{n+\frac{1}{2}}(\vect{x})d\vect{x}}\right) = \Delta t
  \vect{{F}}_{\mathcal{F}}\cdot\vect{V}^{n+\frac{1}{2}}_{\mathcal{F}}
\end{align*}

Finally, using the expression of forces
$\vect{F}_{\mathcal{F}}$, we obtain:
\begin{align*}
  \Delta \mathcal{P}^x_{\mathcal{F}} &= -\Delta
  t\int_{\mathcal{F}}{\overline{p}_xn_{\mathcal{F}}^x}  \\
  \Delta  \mathcal{P}^y_{\mathcal{F}} &= -\Delta t
  \int_{\mathcal{F}}{\overline{p}_yn_{\mathcal{F}}^y} \\
  \Delta \mathcal{E}_{\mathcal{F}} &= -\Delta t
  \vect{V}^{n+\frac{1}{2}}_{\mathcal{F}} \cdot\int_{\mathcal{F}}{\left(
    \begin{array}{c}
      \overline{p}_xn_{\mathcal{F}}^x \\ \overline{p}_yn_{\mathcal{F}}^y
    \end{array}
    \right)}
\end{align*}

Comparing with section \ref{sec:fluide}, the balance of momentum and energy in the fluid domain results in:
\begin{align*}
  \int_{\Omega_f^{n+1}}{\rho^{n+1}\vect{u}^{n+1}}
  +\sum_{\mathcal{F}}{\Delta\vect{\mathcal{P}}_{\mathcal{F}}} &=
  \int_{\Omega_f^{n}}{\rho^{n}\vect{u}^n} \\
  \int_{\Omega_f^{n+1}}{\rho^{n+1}E^{n+1}}
  +\sum_{\mathcal{F}}{\Delta\mathcal{E}_{\mathcal{F}}} &=
  \int_{\Omega_f^{n}}{\rho^{n}E^n}
\end{align*}
This demonstrates the conservation of momentum
and energy for the coupled system.

\section{Conservation of constant flows}
\label{sec:constant}

In this section we analyze the consistency of the coupling method,
in the sense defined in \cite{IBM30}, meaning exact conservation of
uniform flows by the coupling algorithm. Two cases are analyzed. The
first one, also considered in \cite{IBM30}, consists of a solid
immersed in a fluid and moving at the same velocity. This property is
called ``consistency'' in \cite{IBM30}. The second one,
not considered before, demonstrates the correct representation of
the slip boundary condition along walls. These simple cases have
been a guide to design the algorithm, as the preservation of such
flows is a basic criterion for the quality of
the method. \\
In the whole section, we consider a constant fluid state:
 $\rho^n=\rho_0$, $u^n=u_0$, $v^n=v_0$ and $p^n=p_0$ everywhere. The fluxes $f$
 are such that $f_r = f_l = (\rho_0u_0, \rho_0u_0^2+p_0,
 \rho_0u_0v_0, (\rho_0e_0+p_0)u_0)^T$ and $f_t = f_b = (\rho_0v_0,
 \rho_0u_0v_0, \rho_0v_0^2+p_0, (\rho_0e_0+p_0)u_0)^T$.
 In this case, the efficient pressures on the boundary of the solid
are $\overline{p}_x = \overline{p}_y = p_0$.

\subsection{Steady constant flow with moving boundaries}

We consider an arbitrarily shaped rigid body, moving at constant
velocity with no rotation, immersed in a uniform fluid flowing at
the same velocity.

The solid is a closed set, and we denote by $\Omega_s^n$ the solid domain
at initial time. We have:
\begin{equation*}
  \sum_{\mathcal{F}}{S_{\mathcal{F}}\vect{n}_{\mathcal{F}}} = \oint_{\partial\Omega_s^n}{\vect{n}dS} = \vect{0}
\end{equation*}
Using (\ref{eqn:forcepressionx}) and (\ref{eqn:forcepressiony}), we obtain:
\begin{equation*}
  \sum_{\mathcal{F}}{\vect{F}_{\mathcal{F}}} = -\sum_{\mathcal{F}}{p_0 S_{\mathcal{F}}\vect{n}_{\mathcal{F}}} = \vect{0}
\end{equation*}
This induces:
\begin{equation*}
  \vect{V}_i^{n+1} = \vect{V}_i^{n+\frac{1}{2}} = \vect{V}_i^{n} =
  \transpose{(u_0,v_0)}, \;
  \vect{X}_i^{n+1} = \vect{X}_i^{n} + \Delta t\transpose{(u_0,v_0)} 
\end{equation*}
In the same way,
\begin{equation*}
  \sum_{\mathcal{F}}{\vect{\mathcal{M}}_{\mathcal{F}}} =
  -\sum_{\mathcal{F}}{p_0
    S_{\mathcal{F}}\vect{n}_{\mathcal{F}}\wedge(\vect{X}_i^n-\vect{X}_{\mathcal{F}})}
  = p_0\oint_{\partial\Omega_s}{(\vect{X}_i^n-\vect{X})\wedge\vect{n}dS} = \vect{0}
\end{equation*}
The volume swept by boundary $\mathcal{F}$ is $\Delta
tS_{\mathcal{F}}\transpose{(u_0,v_0)}\cdot\vect{n}_{\mathcal{F}}$. Since
the initial state is constant, $\Delta w_{\mathcal{F}}^n$ is given by:
\begin{equation*}
  \Delta w_{\mathcal{F}}^n = \frac{\Delta tS_{\mathcal{F}}}{\Delta
  x\Delta y}(u_0n_{\mathcal{F}}^x+v_0n_{\mathcal{F}}^y)w_0
\end{equation*}
In addition, as the solid translates without
rotation, the normal vector $\vect{n}_{\mathcal{F}}$ to a boundary
$\mathcal{F}$ is constant in time. Using this property in equations
(\ref{eqn:compatibilitex}) and (\ref{eqn:compatibilitey}), we easily
conclude that $(1-\alpha_{\mathcal{C}})\Delta w_{\mathcal{C}}=0$. Thus $w^{n+1}=w^n$, showing that the constant flow is left unchanged by the coupling algorithm and the mixing of small cells.

\subsection{Free slip along a straight boundary}
\label{sec:slip}

We consider an undeformable, fixed solid consisting in a
semi-infinite half-space. The solid boundary is a straight planar
boundary with a constant normal vector $\vect{n}$ such that:
\begin{equation}
  \vect{u}_0\cdot\vect{n} = 0 \label{eqn:normale}
\end{equation}
This initial state describes the free slip of the fluid along the straight
boundary. In the inviscid case, no boundary layer should develop in the
vicinity of the boundary. The conservation of such flows ensures that the boundary is not seen by the fluid as being
artificially rough.

As the solid is fixed, $\alpha_{\mathcal{C}}$ and
$\kappa_{\mathcal{C}}$ remain constant over time and $\Delta
w^{n}_{\mathcal{F}}$ is equal to zero. From equation
(\ref{eqn:fluxmodif2}), and using (\ref{eqn:compatibilitex}),
(\ref{eqn:compatibilitey}) and (\ref{eqn:normale}), the components
of $\Delta w_{\mathcal{C}}$ are calculated as:
\begin{align*}
    (1-\alpha_{\mathcal{C}})\Delta \rho_{\mathcal{C}} &= -\Delta
    t\sum_{\mathcal{F}\in\mathcal{C}}{\frac{S_{\mathcal{F}}}{\Delta
        x\Delta y}\vect{n}\cdot\vect{u}_0} = 0 \\
    (1-\alpha_{\mathcal{C}})\Delta (\rho u)_{\mathcal{C}} &= -\Delta
    t\sum_{\mathcal{F}\in\mathcal{C}}{\frac{S_{\mathcal{F}}}{\Delta
        x\Delta y}\left( (\vect{n}\cdot\vect{u}_0)u_0+p_0n_x \right)} +
    \sum_{\mathcal{F}\in\mathcal{C}}{\frac{\Delta
        tS_{\mathcal{F}}}{\Delta x\Delta y}p_0n_x} = 0 \\
    (1-\alpha_{\mathcal{C}})\Delta (\rho v)_{\mathcal{C}} &= -\Delta
    t\sum_{\mathcal{F}\in\mathcal{C}}{\frac{S_{\mathcal{F}}}{\Delta
        x\Delta y}\left( (\vect{n}\cdot\vect{u}_0)v_0+p_0n_y \right)} +
    \sum_{\mathcal{F}\in\mathcal{C}}{\frac{\Delta
        tS_{\mathcal{F}}}{\Delta x\Delta y}p_0n_y} = 0 \\
  (1-\alpha_{\mathcal{C}})\Delta (\rho E)_{\mathcal{C}} &= -\Delta
    t\sum_{\mathcal{F}\in\mathcal{C}}{\frac{S_{\mathcal{F}}}{\Delta
        x\Delta y}(\vect{n}\cdot\vect{u}_0)(\rho_0e_0+p_0)} = 0
\end{align*}

This shows that the constant flow is preserved by step (5) of the
algorithm. This result is not modified by the mixing procedure. We
thus have shown the exact preservation of the free slip of the fluid
along a straight boundary.

\section{Numerical examples}
\label{sec:numeric}

In the following, we consider a perfect gas, with $\gamma=1.4$. In
all computations the CFL number was fixed equal to 0.5.

\subsection{One-dimensional results}

\begin{figure}[!ht]
\centering \resizebox{0.7\textwidth}{!}{
  \begin{tikzpicture}
    \draw (0,0) node[anchor=south west]{
      \input{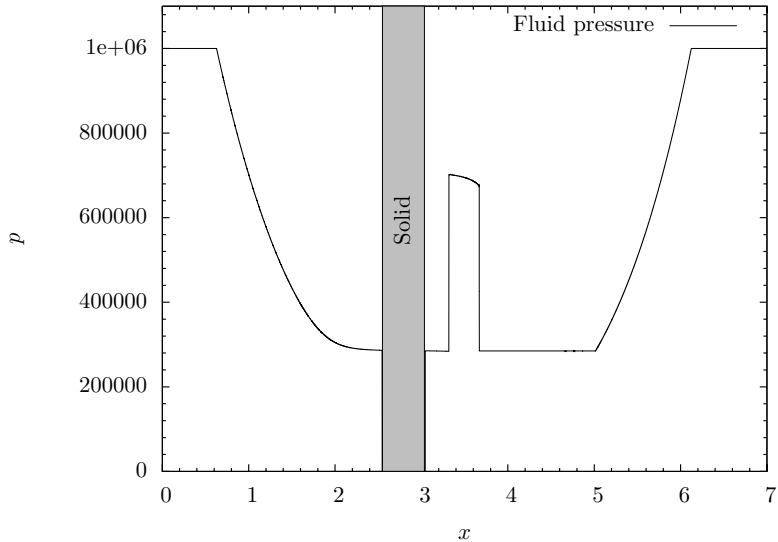}
    };
    \draw[fill=lightgray] (6.21,1.29) -- (6.21,8.55) -- (6.86,8.55) --
    (6.86,1.29) -- cycle;
    \draw (6.76,5.21) node[anchor=south,rotate=90]{Solid};
  \end{tikzpicture}
} \caption{Pressure at time $t=0.003\text{ s}$}\label{fig:pressure}
\end{figure}

\begin{figure}[!ht]
\centering \resizebox{0.7\textwidth}{!}{
  \input{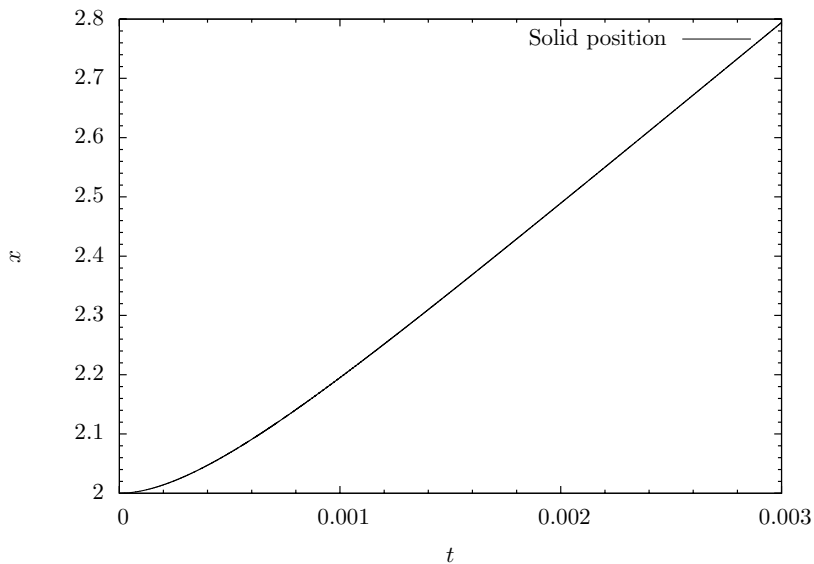}
} \caption{Time evolution of the solid position}\label{fig:position}
\end{figure}

\begin{figure}[!ht]
\centering \resizebox{0.7\textwidth}{!}{
  \includegraphics{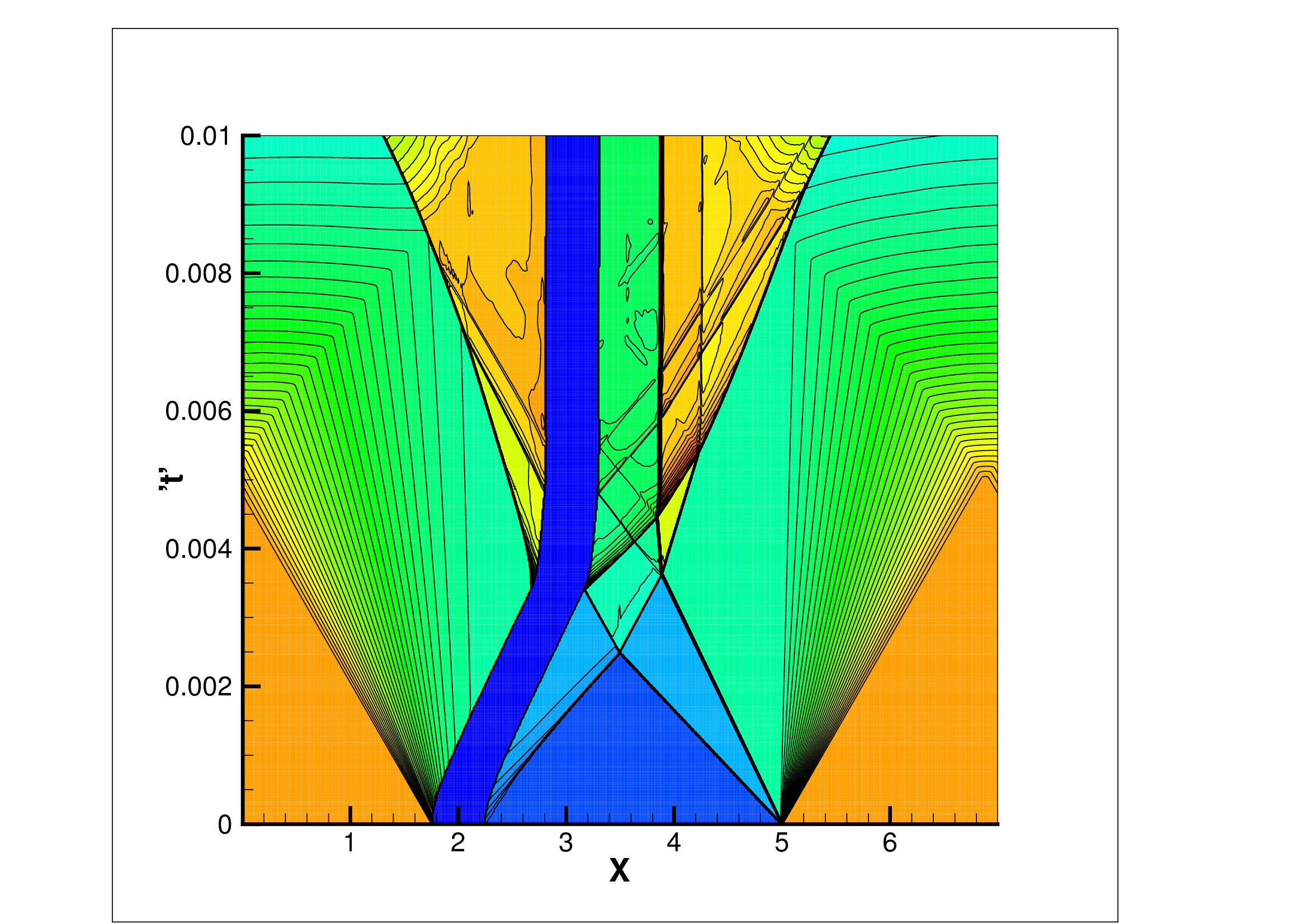}
} \caption{$x-t$ diagram (the position of the solid is in deep
blue)}\label{fig:xt}
\end{figure}

A piston of density $2\text{ kg.m}^{-3}$ and length $0.5\text{ m}$
is initially centered at $x=2\text{ m}$, in a one-dimensional,
$7\text{m}$-long tube, whose ends are connected by periodic boundary
conditions which allow an easier comparison with ALE results. The
gas initial pressure and density are equal to $10^6\text{ Pa}$ and
$10\text{ kg.m}^{-3}$ for $x\leq 2\text{m}$ and $x\geq 5\text{m}$
and to $10^5\text{ Pa}$ and $1\text{ kg.m}^{-3}$ elsewhere. The
system is initially at rest. The initial pressure difference between
the two sides of the piston triggers its movement and the
propagation of waves in the fluid regions (a rarefaction in the left
region and a shock wave in the right one). Wave interactions then
occur at later time. The fluid pressure at time $t=0.003\text{s}$ is
shown in Figure \ref{fig:pressure}, and the trajectory of the solid
is presented in Figure \ref{fig:position}. The $x-t$ diagram over a
longer time (0.01 s) is shown in Figure \ref{fig:xt}.

\begin{figure}[!ht]
\centering \resizebox{0.7\textwidth}{!}{
  \begin{tikzpicture}
    \draw (0,0) node[anchor=south west]{
      \input{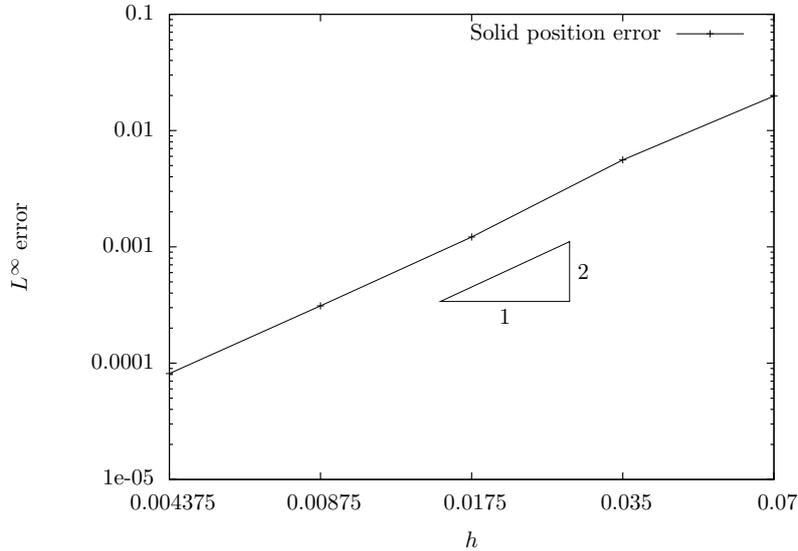}
    };
    \draw (7.,4.0662851) -- (9.,5.) -- (9.,4.0662851) -- cycle;
    \draw (8.,4.0662851) node[anchor=north]{1};
    \draw (9.,4.53) node[anchor=west]{2};
  \end{tikzpicture}
} \caption{Convergence of the solid position
$L^\infty$-error}\label{fig:solid_pos}
\end{figure}

\begin{figure}[!ht]
\centering \resizebox{0.7\textwidth}{!}{
  \begin{tikzpicture}
    \draw (0,0) node[anchor=south west]{
      \input{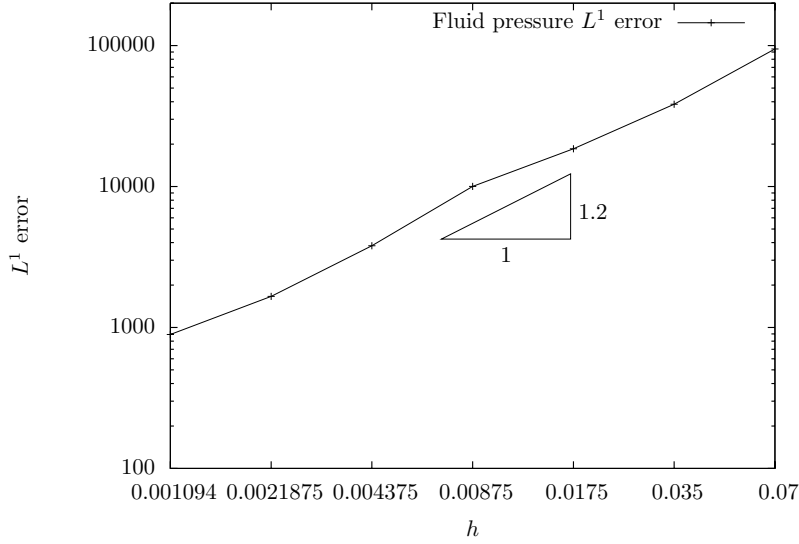}
    };
    \draw (9.,5.885157205) -- (7.,4.8671656) -- (9.,4.8671656) -- cycle;
    \draw (8.,4.8671656) node[anchor=north]{1};
    \draw (9.,5.3) node[anchor=west]{1.2};
  \end{tikzpicture}
} \caption{Convergence of the fluid pressure
$L^1$-error}\label{fig:error_press}
\end{figure}

An ALE computation was done for comparison, using a uniform grid
moving at the solid velocity. The solid position and velocity are
updated using the same second-order Verlet scheme. We compared the
numerical results obtained through the Embedded Boundary method on
100, 200, 400, 800, 1600, 3200, 6400 and 12800  points grids with a
51200-points ALE grid, considered as the reference solution. We
observe a second-order convergence of the solid position (Figure
\ref{fig:solid_pos}) and a super-linear convergence of order 1.2 of
the fluid pressure (Figure \ref{fig:error_press}). The convergence
rate is optimal for the solid (Verlet scheme is second-order
accurate). The convergence rate for the fluid pressure is not
optimal, due to the presence of discontinuities, but is not affected
by the solid coupling.

\subsection{Double Mach reflection}

A Mach 10 planar shock wave reflects on a fixed $30^{\circ}$ wedge,
creating a Mach front, a reflected shock wave, and a contact
discontinuity which develops into a jet along the solid boundary.
This benchmark was first simulated on a Cartesian grid aligned with
the solid boundary, using different finite volume methods
\cite{IBM33,IBM34,Daru2004}. Non-aligned grid methods were also
tested on this benchmark, using Embedded Boundary methods
\cite{IBM18,IBM21}, non-conservative Immersed Boundary methods
\cite{ForrerJeltsch}, $h$-box methods \cite{Helzel}, and kinetic
schemes \cite{KeenKarni}. The position of the tip of the jet is an
important characteristic of the accuracy of the results. The
comparison with grid-aligned results shows that it is better
recovered by conservative methods than non-conservative methods
\cite{IBM18,IBM21}.

\begin{figure}[!ht]
  \centering
  \resizebox{\textwidth}{!}{
    \begin{tikzpicture}
      \draw (0,0) node[anchor=south west]{
       \includegraphics[width=15cm]{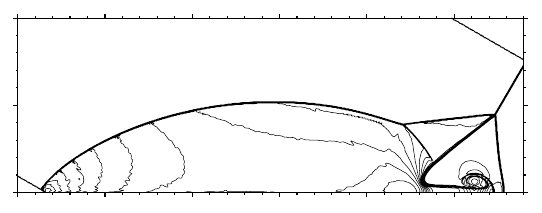}
      };
      \draw (0.65cm,0.5cm) node[anchor=north]{\small 0};
      \draw (3.06cm,0.5cm) node[anchor=north]{\small 0.5};
      \draw (5.47cm,0.5cm) node[anchor=north]{\small 1};
      \draw (7.88cm,0.5cm) node[anchor=north]{\small 1.5};
      \draw (10.29cm,0.5cm) node[anchor=north]{\small 2};
      \draw (12.7cm,0.5cm) node[anchor=north]{\small 2.5};
      \draw (7.64cm,0.1cm) node[anchor=north]{\small x};
      \draw (0.5cm,0.7cm) node[anchor=east]{\small 0};
      \draw (0.5cm,3.11cm) node[anchor=east]{\small 0.5};
      \draw (0.5cm,5.52cm) node[anchor=east]{\small 1};
      \draw (-0.1cm,3.01cm) node[anchor=east]{\rotatebox{90}{\small y}};
    \end{tikzpicture}
  }
  \caption{Aligned case: 30 contours of fluid density from 1.73 to 21, $\varDelta x = \varDelta y = 1/220$} \label{fig:rampe_alignee}
\end{figure}

\begin{figure}[!ht]
  \centering
  \resizebox{\textwidth}{!}{
    \begin{tikzpicture}
      \draw (0,0) node[anchor=south west]{
       \includegraphics[width=15cm]{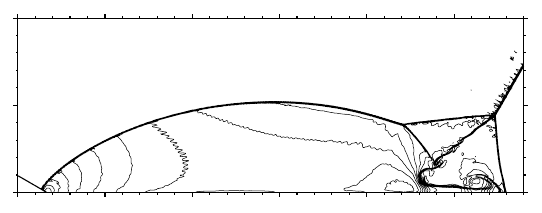}
      };
      \draw (0.65cm,0.5cm) node[anchor=north]{\small 0};
      \draw (3.06cm,0.5cm) node[anchor=north]{\small 0.5};
      \draw (5.47cm,0.5cm) node[anchor=north]{\small 1};
      \draw (7.88cm,0.5cm) node[anchor=north]{\small 1.5};
      \draw (10.29cm,0.5cm) node[anchor=north]{\small 2};
      \draw (12.7cm,0.5cm) node[anchor=north]{\small 2.5};
      \draw (7.64cm,0.1cm) node[anchor=north]{\small x};
      \draw (0.5cm,0.7cm) node[anchor=east]{\small 0};
      \draw (0.5cm,3.11cm) node[anchor=east]{\small 0.5};
      \draw (0.5cm,5.52cm) node[anchor=east]{\small 1};
      \draw (-0.1cm,3.01cm) node[anchor=east]{\rotatebox{90}{\small y}};
    \end{tikzpicture}
  }
  \caption{Non-aligned grid case: 30 contours of fluid density from 1.73 to 21, $\varDelta x = \varDelta y = 1/220$} \label{fig:rampe_inclinee}
\end{figure}

We have simulated the problem on a grid aligned with the wedge
(aligned case, Figure \ref{fig:rampe_alignee}) and on a grid aligned
with the incident shock wave (non-aligned case, Figure
\ref{fig:rampe_inclinee}). The two results are very similar, and
agree with \cite{IBM21,IBM18,IBM33,Daru2004}. One can remark that
all the features of the flow are captured at the correct position in
the non-aligned case. The jet propagates along the wall without
numerical friction due to the conservation of free slip along a
straight boundary (section \ref{sec:slip}). In the principal Mach
stem, the discontinuities are slightly more oscillatory than in the
aligned case. This can be identified as a post-shock oscillation
phenomenon to which Roe's scheme is especially prone (see, for
instance, \cite{Liu96,IBM14}), and is not related to the coupling
method. Nevertheless, the perturbations stay localized in the
vicinity of discontinuities.

\subsection{Lift-off of a cylinder}

\begin{figure}[!ht]
  \centering
  \resizebox{\textwidth}{!}{
    \begin{tikzpicture}
      \draw (0,0) node[anchor=south west]{
       \includegraphics[width=15cm]{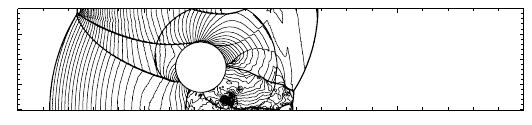}
      };
      \draw (0.6cm,0.5cm) node[anchor=north]{\small 0};
      \draw (4.175cm,0.5cm) node[anchor=north]{\small 0.25};
      \draw (7.75cm,0.5cm) node[anchor=north]{\small 0.5};
      \draw (11.325cm,0.5cm) node[anchor=north]{\small 0.75};
      \draw (14.9cm,0.5cm) node[anchor=north]{\small 1};
      \draw (7.75cm,0.1cm) node[anchor=north]{\small x};
      \draw (0.6cm,0.5cm) node[anchor=east]{\small 0};
      \draw (0.6cm,1.215cm) node[anchor=east]{\small 0.05};
      \draw (0.6cm,1.93cm) node[anchor=east]{\small 0.1};
      \draw (0.6cm,2.645cm) node[anchor=east]{\small 0.15};
      \draw (0.6cm,3.36cm) node[anchor=east]{\small 0.2};
      \draw (0.1cm,1.93cm) node[anchor=east]{\rotatebox{90}{\small y}};
      \draw (7.75cm,3.36cm) node[anchor=south]{\small \bf t=0.14s};
    \end{tikzpicture}
  }
  \resizebox{\textwidth}{!}{
    \begin{tikzpicture}
      \draw (0,0) node[anchor=south west]{
       \includegraphics[width=15cm]{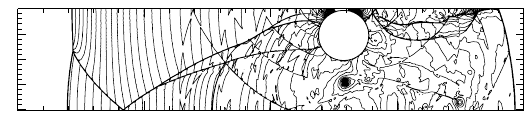}
      };
      \draw (0.6cm,0.5cm) node[anchor=north]{\small 0};
      \draw (4.175cm,0.5cm) node[anchor=north]{\small 0.25};
      \draw (7.75cm,0.5cm) node[anchor=north]{\small 0.5};
      \draw (11.325cm,0.5cm) node[anchor=north]{\small 0.75};
      \draw (14.9cm,0.5cm) node[anchor=north]{\small 1};
      \draw (7.75cm,0.1cm) node[anchor=north]{\small x};
      \draw (0.6cm,0.5cm) node[anchor=east]{\small 0};
      \draw (0.6cm,1.215cm) node[anchor=east]{\small 0.05};
      \draw (0.6cm,1.93cm) node[anchor=east]{\small 0.1};
      \draw (0.6cm,2.645cm) node[anchor=east]{\small 0.15};
      \draw (0.6cm,3.36cm) node[anchor=east]{\small 0.2};
      \draw (0.1cm,1.93cm) node[anchor=east]{\rotatebox{90}{\small y}};
      \draw (7.75cm,3.36cm) node[anchor=south]{\small \bf t=0.255s};
    \end{tikzpicture}
  }
  \caption{60 contours of fluid pressure from 0 to 28 at different
    times, $\varDelta x = \varDelta y = 6.25\times10^{-4}$} \label{fig:lift-off}
\end{figure}

This moving body test case was first proposed in \cite{IBM30}, using
a conservative method. A rigid cylinder of density $7.6
\text{kg.m}^{-3}$ and diameter $0.1\text{ m}$, initially resting on
the lower wall of a $1\text{m}\times0.2\text{m}$ two-dimensional
channel filled with air at standard conditions ($\rho=1\text{
kg.m}^{-3}$, $p=1\text{ Pa}$), is driven and lifted upwards by a
Mach 3 shock wave. Gravity is not taken into account. The problem
was simulated in \cite{ForrerBerger,IBM37,IBM10,IBM22}. In Figure
\ref{fig:lift-off}, we present our results on a uniform $1600\times
320$ grid at times $0.14\text{ s}$ and $0.255\text{ s}$. The
cylinder was approximated by a polygon with 1240 faces.

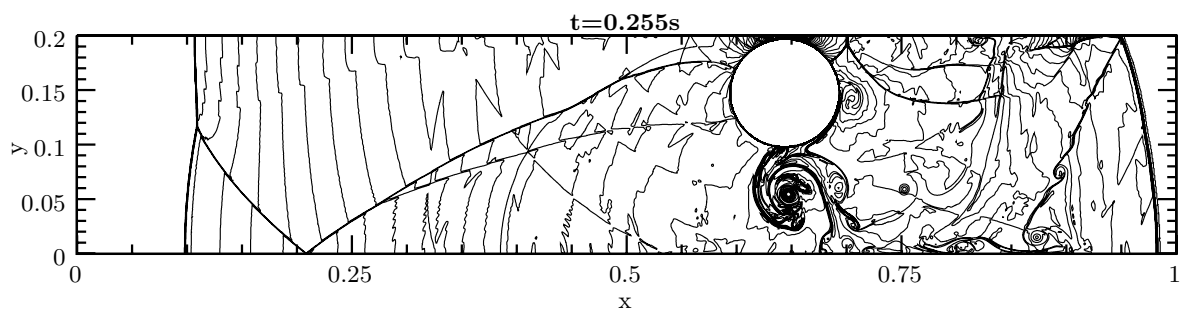
\begin{figure}[!ht]
  \centering
  \resizebox{\textwidth}{!}{
    \begin{tikzpicture}
      \draw (0,0) node[anchor=south west]{
       \includegraphics[width=15cm]{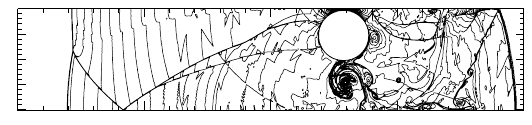}
      };
      \draw (0.6cm,0.5cm) node[anchor=north]{\small 0};
      \draw (4.175cm,0.5cm) node[anchor=north]{\small 0.25};
      \draw (7.75cm,0.5cm) node[anchor=north]{\small 0.5};
      \draw (11.325cm,0.5cm) node[anchor=north]{\small 0.75};
      \draw (14.9cm,0.5cm) node[anchor=north]{\small 1};
      \draw (7.75cm,0.1cm) node[anchor=north]{\small x};
      \draw (0.6cm,0.5cm) node[anchor=east]{\small 0};
      \draw (0.6cm,1.215cm) node[anchor=east]{\small 0.05};
      \draw (0.6cm,1.93cm) node[anchor=east]{\small 0.1};
      \draw (0.6cm,2.645cm) node[anchor=east]{\small 0.15};
      \draw (0.6cm,3.36cm) node[anchor=east]{\small 0.2};
      \draw (0.1cm,1.93cm) node[anchor=east]{\rotatebox{90}{\small y}};
      \draw (7.75cm,3.36cm) node[anchor=south]{\small \bf t=0.255s};
    \end{tikzpicture}
  }
  \caption{60 contours of fluid density from 0 to 12 at final time, $\varDelta x = \varDelta y = 6.25\times10^{-4}$} \label{fig:densite_lift-off}
\end{figure}

Our results agree well on the position of the solid and of the
shocks with Arienti et al. \cite{IBM37} and Hu et al. \cite{IBM10}.
However, some differences should be noted. First of all, some
reflected shock waves in our results seem to lag slightly behind
their position in previous results. This difference might be caused
by small differences in the final position of the solid. Hu et al.
\cite{IBM10} also discuss the presence of a strong vortex under the
cylinder in the results of Forrer and Berger \cite{ForrerBerger}.
They dismiss it as an effect of the space-time splitting scheme
employed which affects the numerical dissipation. We also obtain
this vortex, which does not disappear as we refine the mesh. We
rather believe that this vortex is associated with a
Kelvin-Helmholtz instability of the contact discontinuity present
under the cylinder (Figure \ref{fig:densite_lift-off}).

In Figures \ref{fig:position_cylindre_x} and
\ref{fig:position_cylindre_y} we present convergence results on the
final position of the center of mass of the cylinder, compared to
those of Hu et al. \cite{IBM10}. We observe that our results exhibit
a fast convergence process,  which is not the case in \cite{IBM10}.
Let us note that no exact solution exists for the final position of
the cylinder. The final position we found is however in the same
range as in \cite{IBM10}. The results also compare well with Arienti
et al. \cite{IBM37}. The improvement lies in the combination of the
conservative interface method \cite{IBM10} with a conservative
coupling and a second-order time-scheme for the rigid body motion.
For this difficult case, the maximal conservation relative errors
due to coupling were bounded by $4.10^{-6}$ over the whole
simulation time, and no drift was observed.

\begin{figure}[!ht]
  \centering
  \resizebox{0.8\textwidth}{!}{
    \begin{tikzpicture}
      \draw (0,0) node[anchor=south west]{
        \input{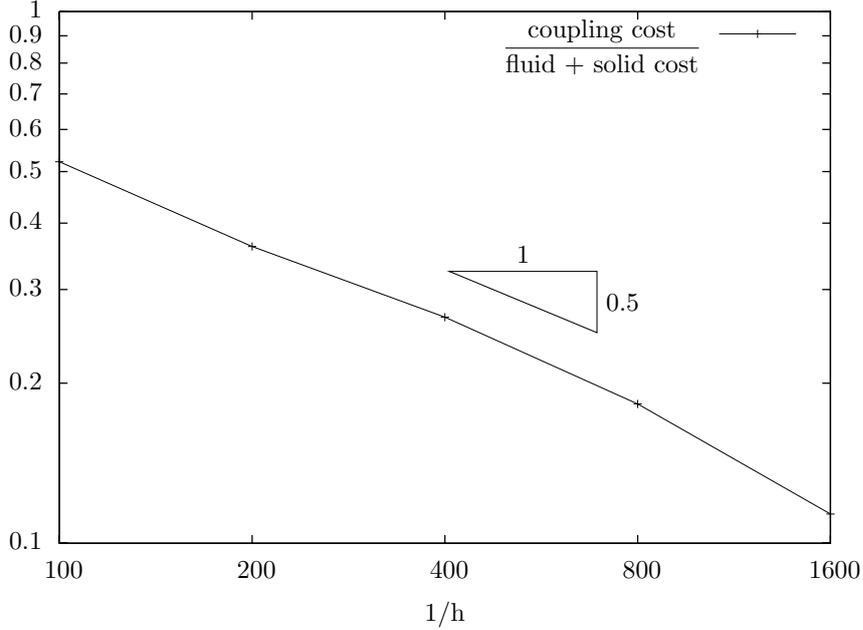}
      };
      \draw (7.,5.) -- (9.,5.) -- (9.,4.163403) -- cycle;
      \draw (8.,5.) node[anchor=south]{1};
      \draw (9.,4.5817) node[anchor=west]{0.5};
    \end{tikzpicture}
  }
  \caption{Ratio between the coupling computation cost and the fluid
    and solid costs} \label{fig:temps_calcul_liftup}
\end{figure}

In Figure \ref{fig:temps_calcul_liftup} we present the relative
computational cost of the coupling. The relative cost is defined as
the ratio of the computational times dedicated to the coupling
method and to the fluid and solid methods. In the rigid body case,
the cost of the solid method is very low compared to that of the
fluid method. As the coupling method is explicit and local, the
computational cost is located on a manifold  one dimension lower
than the dimension of the whole space. In the two-dimensional case,
the coupling is on a one-dimensional manifold. Indeed, we observe in
Figure \ref{fig:temps_calcul_liftup} that the relative cost of the
coupling decreases as the grid is refined, with a slope of $0.5$,
and that the coupling cost remains lower than the fluid and solid
costs, amounting to approximately 10--20\% for the grids yielding
sufficient accuracy.

\begin{figure}[!ht]
  \centering
  \resizebox{0.8\textwidth}{!}{
    \input{convergencex.tex}
  }
  \caption{Comparison of the convergence of the horizontal position of the center
    of mass of the cylinder with Hu et al. \cite{IBM10}} \label{fig:position_cylindre_x}
\end{figure}

\begin{figure}[!ht]
  \centering
  \resizebox{0.8\textwidth}{!}{
    \input{convergencey.tex}
  }
  \caption{Comparison of the convergence of the vertical position of the center
    of mass of the cylinder with Hu et al. \cite{IBM10}} \label{fig:position_cylindre_y}
\end{figure}

\subsection{Flapping doors}

\begin{figure}[!ht]
  \centering
  \resizebox{1.\textwidth}{!}{
    \begin{tikzpicture}
      \draw (0,0) node[anchor=south west]{
       \includegraphics[width=15cm]{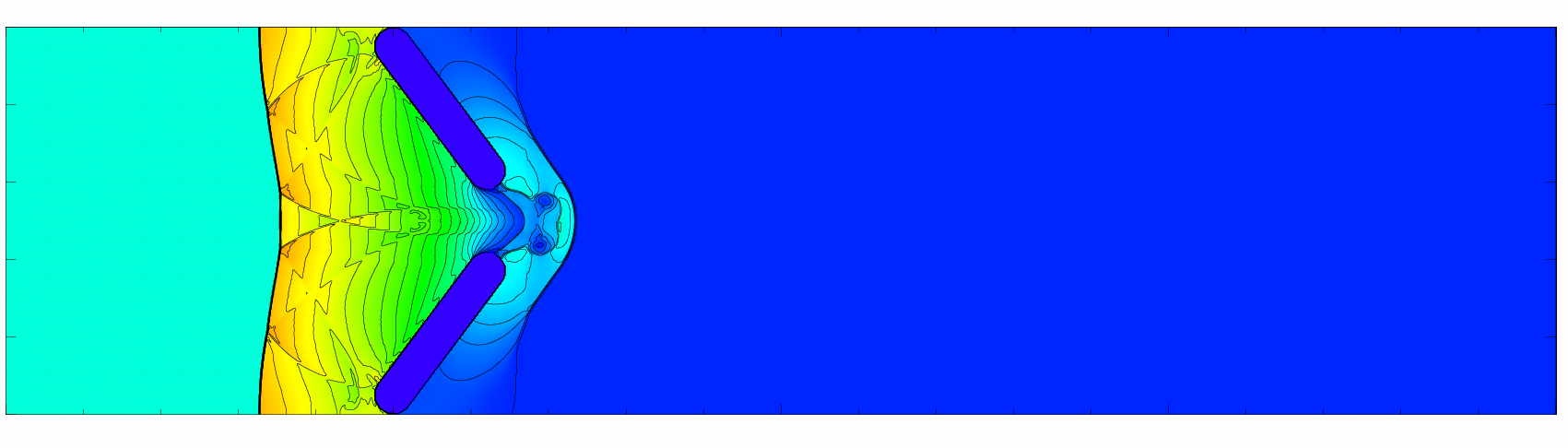}
      };
      \draw (0.2cm,0.5cm) node[anchor=north]{\small 0};
      \draw (3.9cm,0.5cm) node[anchor=north]{\small 0.5};
      \draw (7.6cm,0.5cm) node[anchor=north]{\small 1};
      \draw (11.3cm,0.5cm) node[anchor=north]{\small 1.5};
      \draw (15.cm,0.5cm) node[anchor=north]{\small 2};
      \draw (7.6cm,0.1cm) node[anchor=north]{\small x};
      \draw (0.2cm,0.5cm) node[anchor=east]{\small 0};
      \draw (0.2cm,1.215cm) node[anchor=east]{\small 0.1};
      \draw (0.2cm,1.93cm) node[anchor=east]{\small 0.2};
      \draw (0.2cm,2.645cm) node[anchor=east]{\small 0.3};
      \draw (0.2cm,3.36cm) node[anchor=east]{\small 0.4};
      \draw (0.2cm,4.075cm) node[anchor=east]{\small 0.5};
      \draw (-0.3cm,2.29cm) node[anchor=east]{\rotatebox{90}{\small y}};
      \draw (15.1cm,2.29cm) node[anchor=west]{\small (a)};
    \end{tikzpicture}
  }
%
  \resizebox{1.\textwidth}{!}{
    \begin{tikzpicture}
      \draw (0,0) node[anchor=south west]{
       \includegraphics[width=15cm]{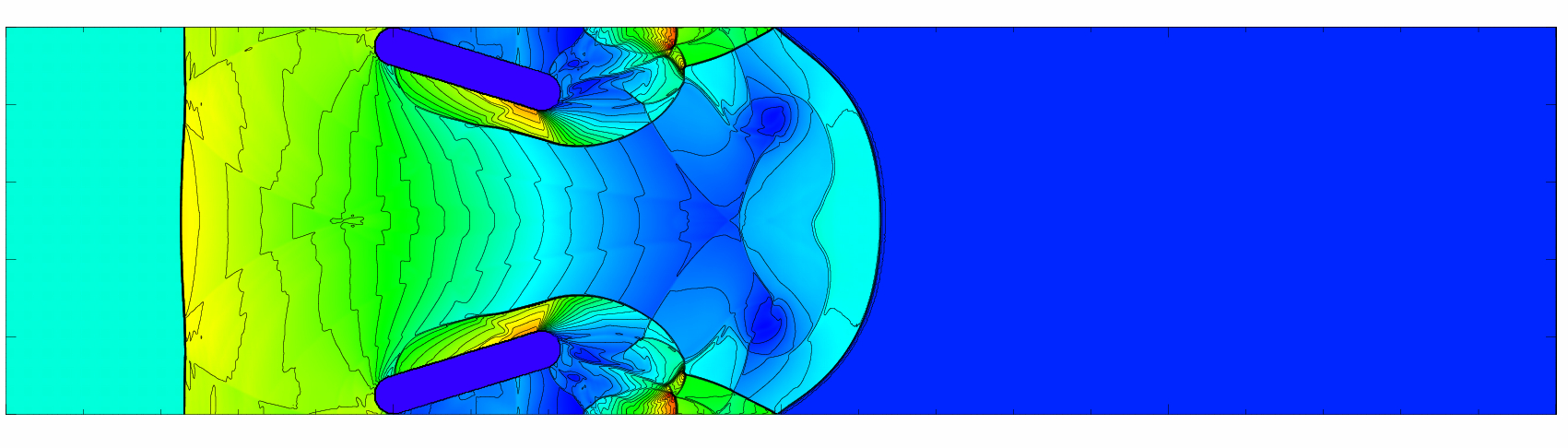}
      };
      \draw (0.2cm,0.5cm) node[anchor=north]{\small 0};
      \draw (3.9cm,0.5cm) node[anchor=north]{\small 0.5};
      \draw (7.6cm,0.5cm) node[anchor=north]{\small 1};
      \draw (11.3cm,0.5cm) node[anchor=north]{\small 1.5};
      \draw (15.cm,0.5cm) node[anchor=north]{\small 2};
      \draw (7.6cm,0.1cm) node[anchor=north]{\small x};
      \draw (0.2cm,0.5cm) node[anchor=east]{\small 0};
      \draw (0.2cm,1.215cm) node[anchor=east]{\small 0.1};
      \draw (0.2cm,1.93cm) node[anchor=east]{\small 0.2};
      \draw (0.2cm,2.645cm) node[anchor=east]{\small 0.3};
      \draw (0.2cm,3.36cm) node[anchor=east]{\small 0.4};
      \draw (0.2cm,4.075cm) node[anchor=east]{\small 0.5};
      \draw (-0.3cm,2.29cm) node[anchor=east]{\rotatebox{90}{\small y}};
      \draw (15.1cm,2.29cm) node[anchor=west]{\small (b)};
    \end{tikzpicture}
  }
%
  \resizebox{1.\textwidth}{!}{
    \begin{tikzpicture}
      \draw (0,0) node[anchor=south west]{
       \includegraphics[width=15cm]{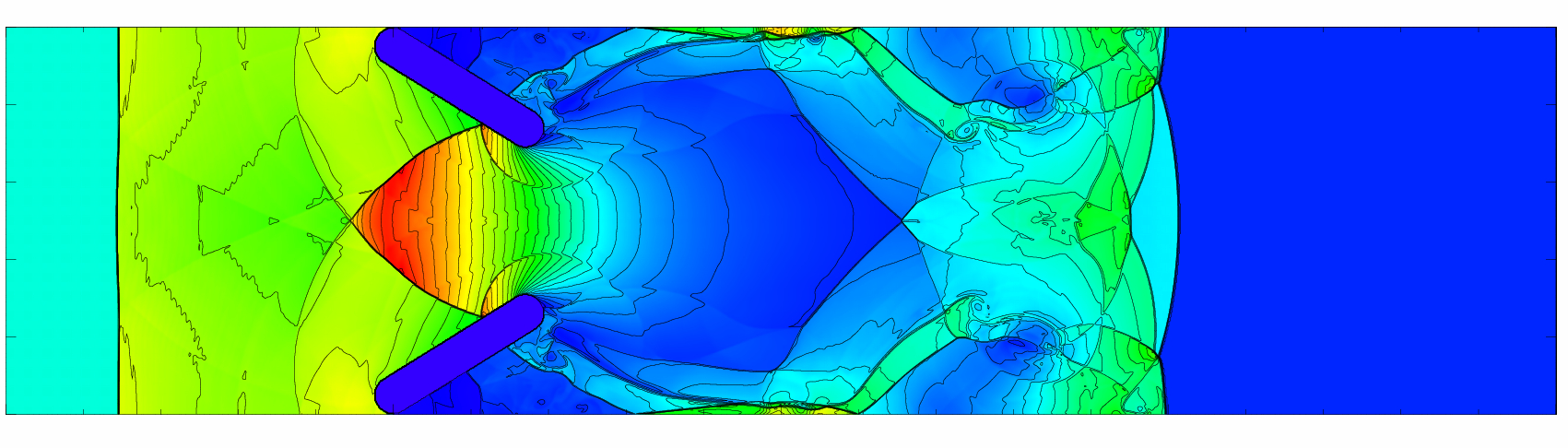}
      };
      \draw (0.2cm,0.5cm) node[anchor=north]{\small 0};
      \draw (3.9cm,0.5cm) node[anchor=north]{\small 0.5};
      \draw (7.6cm,0.5cm) node[anchor=north]{\small 1};
      \draw (11.3cm,0.5cm) node[anchor=north]{\small 1.5};
      \draw (15.cm,0.5cm) node[anchor=north]{\small 2};
      \draw (7.6cm,0.1cm) node[anchor=north]{\small x};
      \draw (0.2cm,0.5cm) node[anchor=east]{\small 0};
      \draw (0.2cm,1.215cm) node[anchor=east]{\small 0.1};
      \draw (0.2cm,1.93cm) node[anchor=east]{\small 0.2};
      \draw (0.2cm,2.645cm) node[anchor=east]{\small 0.3};
      \draw (0.2cm,3.36cm) node[anchor=east]{\small 0.4};
      \draw (0.2cm,4.075cm) node[anchor=east]{\small 0.5};
      \draw (-0.3cm,2.29cm) node[anchor=east]{\rotatebox{90}{\small y}};
      \draw (15.1cm,2.29cm) node[anchor=west]{\small (c)};
    \end{tikzpicture}
  }
%
  \resizebox{1.\textwidth}{!}{
    \begin{tikzpicture}
      \draw (0,0) node[anchor=south west]{
       \includegraphics[width=15cm]{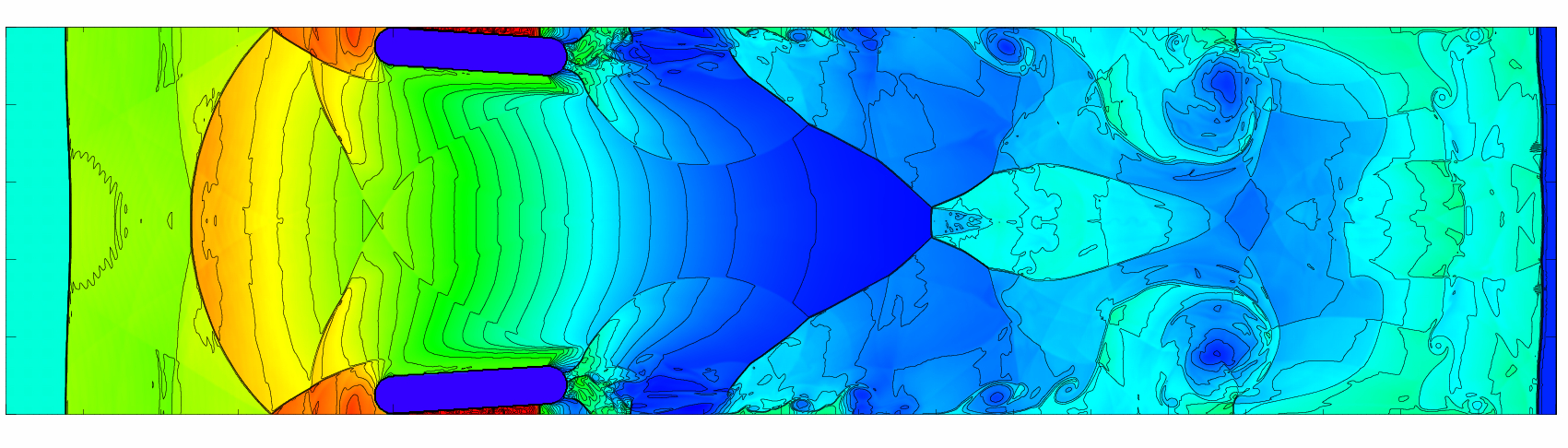}
      };
      \draw (0.2cm,0.5cm) node[anchor=north]{\small 0};
      \draw (3.9cm,0.5cm) node[anchor=north]{\small 0.5};
      \draw (7.6cm,0.5cm) node[anchor=north]{\small 1};
      \draw (11.3cm,0.5cm) node[anchor=north]{\small 1.5};
      \draw (15.cm,0.5cm) node[anchor=north]{\small 2};
      \draw (7.6cm,0.1cm) node[anchor=north]{\small x};
      \draw (0.2cm,0.5cm) node[anchor=east]{\small 0};
      \draw (0.2cm,1.215cm) node[anchor=east]{\small 0.1};
      \draw (0.2cm,1.93cm) node[anchor=east]{\small 0.2};
      \draw (0.2cm,2.645cm) node[anchor=east]{\small 0.3};
      \draw (0.2cm,3.36cm) node[anchor=east]{\small 0.4};
      \draw (0.2cm,4.075cm) node[anchor=east]{\small 0.5};
      \draw (-0.3cm,2.29cm) node[anchor=east]{\rotatebox{90}{\small y}};
      \draw (15.1cm,2.29cm) node[anchor=west]{\small (d)};
    \end{tikzpicture}
  }
  \caption{Density contours at times $t=0.125$s (a), $t=0.25$s (b), $t=0.375$s (c) and $t=0.5$s (d)} \label{fig:flap}
\end{figure}
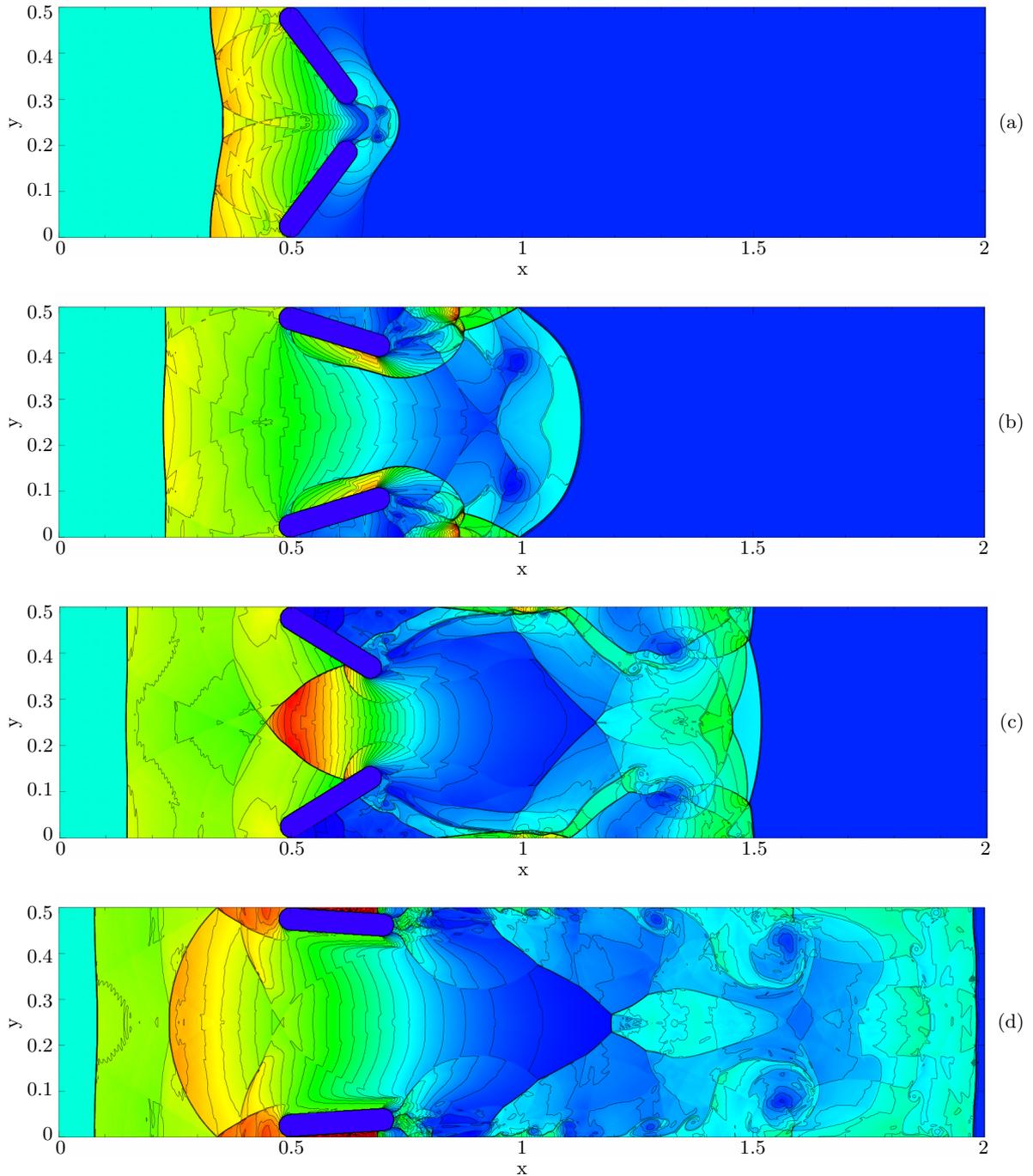

We propose this new fluid-structure interaction case as a
demonstration of the robustness of our approach and as a first step
towards fracture and impact simulations. The flapping doors case
involves separating or closing solid boundaries, with cells
including several moving boundaries.  The algorithm is shown to be
able to deal with such difficulties. Two doors initially close a
canal and are impacted from the left by a Mach 3 shock. The canal
consists of two fixed rigid walls, 2m long and 0.5m apart. Each door
consists of a 0.2-m long and 0.05m-wide rectangle, completed at both
ends by a half-circle of diameter 0.05m. The doors are respectively
fixed on points $(0.5,0.025)$ and $(0.5,0.475)$ at the center of the
half-circles. They can rotate freely around these points. The Mach 3
shock is initially located at $x=0.43$m. The density of the solid is
$0.1$ kg.m$^{-3}$ and the pre- and post-shock state of the fluid are
$(\rho,u,v,p) = (1,0,0,1)$ and
$(\rho,u,v,p)=(3.857,2.6929,0,10.333)$. In Figure \ref{fig:flap} we
show the density field obtained using a 1600x400 grid, at times
$0.125$s, $0.25$s, $0.375$s and $0.5$s. After the incident shock
hits the doors, it reflects to the left and the doors open due to
the high rise in pressure. The opening of the doors produces  a jet
preceded by a shock wave propagating to the right. Then complex
interactions of waves occur, due to door movements and interaction
with the walls. Kelvin-Helmholtz instabilities of contact lines can
be observed at t=0.5s. It is worth noting that symmetry of the flow
about the centerline of the canal axis is remarkably well preserved
by the coupling method.

\begin{figure}[!ht]
  \centering
  \resizebox{0.8\textwidth}{!}{
    \input{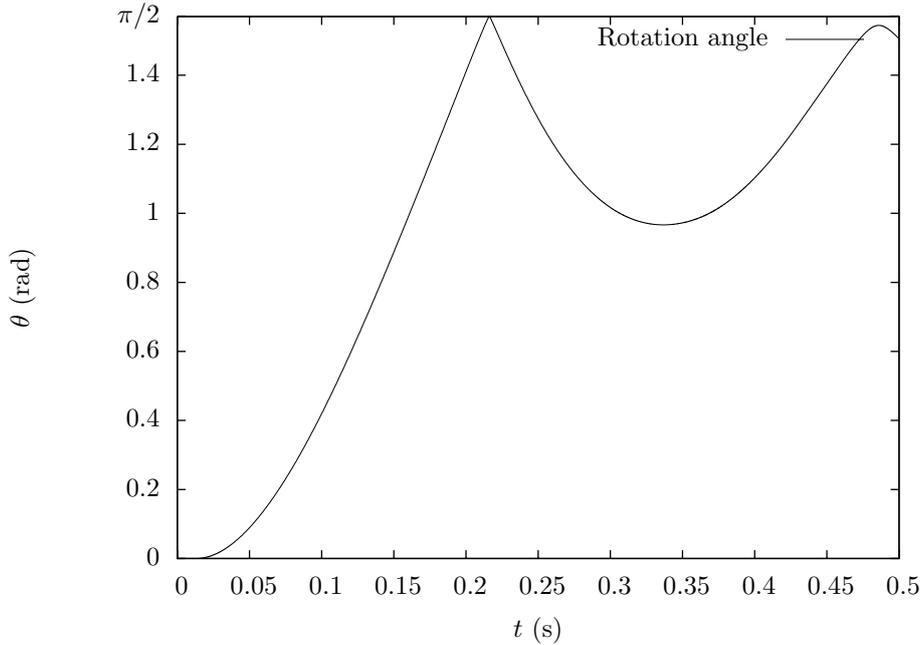}
  }
  \caption{Time evolution of the rotation of a door} \label{fig:rotation_flap}
\end{figure}

As the doors remain tangent to the canal walls during their
rotation, the fluid cannot pass between the wall and a door at its
hinge. When the doors approach the walls at maximum rotation, the
fluid is compressed, and eventually pushes them back. This is
observed at time $t = 0.2162$s and $t = 0.486$s in Figure
\ref{fig:rotation_flap}, which presents the time evolution of the
doors rotation angle. In the first case, the distance between each
door's straight boundary and the wall is less than 0.002m, while the
size of a fluid cell is $\Delta x=0.00125$m. The method is able to
deal with the fact that most cells along the wall are cut by the
moving boundary and contain several moving boundaries. Treating this
test case with an ALE method would require several remeshings in the
course of the simulation, especially in the initial separation of
the tangent door tips, and when the doors approach the walls.

\section{Conclusion}

We have presented a new coupling algorithm between a compressible fluid
flows and a rigid body using an Embedded Boundary method.
This explicit algorithm has the advantage of preserving the usual
CFL stability condition: the time-step can be taken as the minimum
of the full cell size fluid and solid time-steps. The combination of
the Embedded Boundary method for the fictitious fluid domain and of
the coupling strategy ensures the conservation of fluid mass and the
balance of momentum and energy between fluid and solid. In addition,
the exact conservation of the two constant states described in
section \ref{sec:constant} gives good insight on the consistency of
the method: we prove the conservation of a constant flow in which a
solid moves at the same velocity and the fact that
the treatment at the boundary introduces no spurious roughness or
boundary layers. The numerical examples suggest the second-order
convergence of the solid position and the super-linear convergence
of the fluid state in $L^1$ norm, while our results on
two-dimensional benchmarks agree very well with body-fitted methods
and improve Immersed Boundary results. We are also capable of
dealing with solid boundaries moving close to each other, which is
promising for impact simulations. The method is computationally
efficient, as the coupling adds an integration on a space one
dimension smaller than the fluid and solid computation spaces.  The
present method is therefore perfectly liable to be extended to a
deformable solid, and was designed to extend naturally in three
space dimensions. The remaining difficulty is the ability to define
and track the solid boundary surrounding a Discrete Element
assembly.

\FloatBarrier

\bibliographystyle{plain}
\bibliography{jcp2010}

\end{document}